\documentclass[12pt]{article}
\usepackage{amsmath, amssymb, amsthm, amsfonts}
\usepackage{parskip}
\usepackage{geometry}
\usepackage{mathtools}
\usepackage{hyperref}
\usepackage{outlines}
\usepackage{cite}
\usepackage{braket}
\usepackage{float}
\usepackage{caption}
\usepackage{color}

\usepackage[inline]{asymptote}

\title{On graphs with large third eigenvalue}
\author{Giacomo Leonida\thanks{\href{mailto:gl510@cam.ac.uk}{gl510@cam.ac.uk}, \href{mailto:sl2190@cam.ac.uk}{sl2190@cam.ac.uk}. Department of Pure Mathematics and Mathematical Statistics
(DPMMS), University of Cambridge, United Kingdom.}\; and Sida Li\footnotemark[1]}
\date{}

\newcommand{\floor}[1]{\lfloor{#1\rfloor}}
\newcommand{\ceil}[1]{\lceil{#1\rceil}}

\newtheorem{thm}{Theorem}[section]
\newtheorem{dfn}[thm]{Definition}
\newtheorem{lmm}[thm]{Lemma}
\newtheorem{cor}[thm]{Corollary}

\newtheorem{qun}[thm]{Question}
\newtheorem{conj}[thm]{Conjecture}
\newtheorem*{claim}{Claim}

\begin{document}
\maketitle
\begin{abstract}
    Given a graph $G$, let $\lambda_3$ denote the third largest eigenvalue of its adjacency matrix. In this paper, we prove various results towards the conjecture that $\lambda_3(G) \le \frac{|V(G)|}{3}$, motivated by a question of Nikiforov. We generalise the known constructions that yield $\lambda_3(G) = \frac{|V(G)|}{3} - 1$ and prove the inequality holds for $G$ strongly regular, a regular line graph or a Cayley graph on an abelian group. We also consider the extended problem of minimising $\lambda_{n-1}$ on weighted graphs and reduce the existence of a minimiser with simple final eigenvalue to a vertex multiplication of a graph on 11 vertices. We prove that the minimal $\lambda_{n-1}$ over weighted graphs is at most $O(\sqrt{n})$ from the minimal $\lambda_{n-1}$ over unweighted graphs. 
\end{abstract}

\section{Introduction}

For a graph $G$, let $A(G)$ denote its adjacency matrix and $\lambda_k$ the $k$-th largest eigenvalue of $A(G)$. Following Nikiforov's notation in \cite{NikiforovEvals}, define $c_k = \sup\{\lambda_k(G)/n : |V(G)| = n \ge k\}$. Nikiforov proved the following upper bound on $c_k$.

\begin{thm}[Nikiforov \cite{NikiforovEvals}]\label{thm:nikiforov-upper}
    For $k \ge 2$,
    \[c_k \le \frac{1}{2\sqrt{k-1}}.\]
\end{thm}

This gives that $c_2 = \frac{1}{2}$, as the construction $kK_{\frac{n}{k}}$ has $\lambda_k(kK_{\frac{n}{k}}) = \frac{n}{k} - 1$, so $c_k \ge \frac1k$. Nikiforov provided a further construction yielding $c_k \ge \frac{1}{k - \frac12}$ for $k \ge 5$. Additionally, for sufficiently large $k$, Nikiforov used the family of strongly-regular graphs introduced by Taylor in \cite{TaylorGraphs1, TaylorGraphs2}, widely known as Taylor graphs, for a closely matching lower bound to Theorem \ref{thm:nikiforov-upper}.

\begin{thm}[Nikiforov \cite{NikiforovEvals}]
    There exists a $k_0$ such that if $k > k_0$, then
    \[c_k > \frac{1}{2\sqrt{k-1} + \sqrt[3]{k}}.\]
\end{thm}

Out of the two remaining cases $k \in \{3, 4\}$, recently Linz \cite{LinzImproved} tightened the lower bound on $k = 4$ with the closed blow-up of the icosahedral graph.

\begin{thm}[Linz \cite{LinzImproved}]
    \[c_4 \ge \frac{1+\sqrt{5}}{12} \approx 0.26967 > \frac14.\]
\end{thm}

Thus, we have that $c_k > \frac1k$ for $k \ge 4$ and $c_k = \frac1k$ for $k \le 2$. \\

In this paper, we provide evidence for $c_3 = \frac13$, providing a partial affirmative response to Nikiforov's original question.

\begin{qun}[Nikiforov\cite{NikiforovEvals}]\label{qun:c_3}
    Is it true that $c_3 = \frac13$?
\end{qun}

First, we demonstrate a family of graphs, denoted $H_{a,b}$, satisfying $\lambda_3(G) = \frac{|V(G)|}{3} - 1$, which generalises the only previously-known constructions of $3K_n$ and the closed blow-up of $C_6$. In particular, these include the first non-regular examples of such graphs and we propose the following conjecture.

\begin{conj}
    If $\lambda_3(G) = \frac{|V(G)|}{3} - 1$, then $G$ is isomorphic to $H_{a,b}$ for some non-negative integers $a,b$. 
\end{conj}

We then prove that multiple families of graphs, which provide the best known lower bounds for larger values of $k$, have $\lambda_3(G) \le \frac{|V(G)|}{3}$, further suggesting $c_3 = \frac13$. In particular, we prove and generalise the following theorem.

\begin{thm}\label{th:classes}
    If $G$ is strongly regular, a regular line graph, or a Cayley graph on an abelian group, then $\lambda_3(G) \le \frac{|V(G)|}{3}$. 
\end{thm}

Note that within the table of best-known lower bounds compiled by Linz \cite{LinzImproved}, Taylor graphs are strongly regular, Paley graphs are Cayley over $\mathbb{Z}_n$ and both the Petersen graph and $C_6$ are regular line graphs. In addition, we generalise the methods of resolving Cayley graphs on an abelian group to semi-Cayley graphs on an abelian group, as well as vertex-transitive graphs with abelian automorphism group. 

We also consider the extended problem of bounding $\lambda_k$ for weighted graphs with weights in $[0,1]$. This was motivated by noting that weighted graphs can be interpreted as graphons in the sense of Lov\'{a}sz and Szegedy \cite{LovaszGraphon} and thus we can approximate these with a convergent sequence of pseudo-random graphs. 

In particular, we consider the problem of minimising $\lambda_{n-1}$ over such weighted graphs, which by Weyl's inequality (see \cite{SpectraTextbook} p.19), is analogous to the problem of maximising $\lambda_3$. By considering the effects of suitable edge operations, we are able to derive various structural results on the minimisers of $\lambda_{n-1}$. Namely, we prove that the extension of Question \ref{qun:c_3} to weighted graphs does not affect the desired constant in the following sense.
\begin{thm}\label{thm:weighted-equivalent-intro}
    For $n \ge 2$, if $\mathcal{W}_n$ denotes the set of weighted graphs on $n$ vertices with weights in $[0,1]$ and $\mathcal{G}_n$ denotes the set of unweighted graphs on $n$ vertices, then:
\begin{equation*}
    \min_{G \in \mathcal{G}_n} \lambda_{n-1}(G) - 2\sqrt{n} \le \inf_{G\in \mathcal{W}_n} \lambda_{n-1}(G).
\end{equation*}
\end{thm}
Furthermore, we explore the existence of a minimiser of $\lambda_{n-1}$ with simple final eigenvalue. Let $\mathcal{S}_n$ denote the set of $n\times n$ symmetric matrices with entries in $[0,1]$. We prove the following as part of Theorem \ref{thm:weighted-equivalent-intro}. 
\begin{thm} If there exists a minimiser of $\lambda_{n-1}$ over $\mathcal{S}_n$ with simple final eigenvalue, then there exists an $11\times11$ matrix $T$ with entries in $\{0,1\}$ and a vertex multiplication $\tilde{T}$ of $T$ with:
\begin{equation*}
    \lambda_{n-1}(\tilde{T}) - 2\sqrt{n} \le  \inf_{S\in \mathcal{S}_n} \lambda_{n-1}(S).
\end{equation*}
\end{thm}

We conclude by providing some small structural results regarding minimisers of $\lambda_{n-1}$ with repeated final eigenvalue. 
\section{Graphs with $\lambda_3(G) = \frac{|V(G)|}{3}-1$}

We construct a family of graphs which generalises the $3K_m$ and $C_6^{[m]}$ (the closed blowup of the 6-cycle, as described by Linz \cite{LinzImproved}) constructions, with third eigenvalue $\frac{|V(G)|}{3}-1$. The closed blow-ups of these graphs have third eigenvalue proportion converging to $\frac13$.

A closed vertex multiplication of a graph $G$ is constructed by replacing each vertex of $G$ with a clique and each edge of $G$ with a complete bipartite graph between the endpoint cliques. For a graph $G$ with vertices $v_1,\dots, v_n$, the closed vertex multiplication of $G$ by $[a_1,\dots, a_n]$ is the closed vertex multiplication where the vertex $v_i$ is replaced by a clique of size $a_i$. If the cliques are chosen to be the same size, this reduces to a closed blow-up.

\begin{dfn}\label{def:closed-vert}
    Let $H_{a,b}$ denote the closed vertex multiplication of $C_6$ by $[a, b, a, b, a, b]$. 
\end{dfn}
\begin{figure}[ht!]
    \centering
    \includegraphics[scale=0.9]{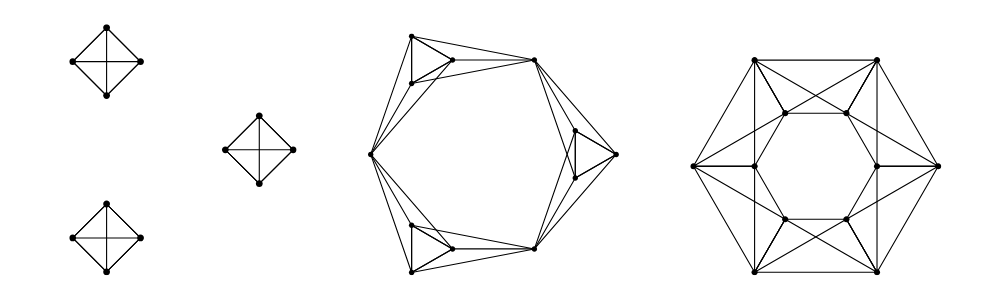}
    \caption{$H_{4,0}$, $H_{3,1}$ and $H_{2,2}$ respectively.}
    \label{fig:C6VertMult}
\end{figure}
In particular, we remark that $H_{n,0}$ is isomorphic to $3K_n$ and $H_{n,n}$ is isomorphic to $C_6^{[n]}$. 
\begin{thm}\label{thm:new-family} We have $\lambda_3(H_{a,b}) = a+b-1 = \frac{|V(H_{a,b})|}{3}-1$. Hence, 
    $\lim_{t \to \infty} \frac{\lambda_3(H_{a,b}^{[t]})}{3(a+b)t} = \frac13$. 
\end{thm}
\begin{proof}
    Add a self-loop to each vertex of $H_{a,b}$. Then, there exists an equitable partition of this graph given by the 6 cliques. The divisor matrix of this equitable partition is:
    \[\begin{pmatrix}
        a & b & 0 & 0 & 0 & b \\
        a & b & a & 0 & 0 & 0 \\
        0 & b & a & b & 0 & 0 \\
        0 & 0 & a & b & a & 0 \\
        0 & 0 & 0 & b & a & b \\
        a & 0 & 0 & 0 & a & b
    \end{pmatrix}.\]
    Hence, $A(H_{a,b}) + I$ has spectrum:
    \[\frac{a+b}{2} + \frac{\sqrt{a^2 + 14ab + b^2}}{2}, a+b, a+b, 0, \dots, 0, \frac{a+b}{2} - \frac{\sqrt{a^2 + 14ab + b^2}}{2},\]
    and $\lambda_3(H_{a,b}) = a+b-1 = \frac{|V(H_{a,b})|}{3}-1$. Thus, $\lim_{t \to \infty} \frac{\lambda_3(H_{a,b}^{[t]})}{3(a+b)t} = \lim_{t \to \infty} \frac{(a+b)t - 1}{3(a+b)t} = \frac13$.
\end{proof} 

\section{Bounds on the third eigenvalue of graphs of given type}
In this section, we explore bounds on the third eigenvalue for various classes of graphs, proving Theorem \ref{th:classes}. We start with strongly regular graphs.
\subsection{Strongly regular graphs}
Several known examples of graphs with large $k$-th eigenvalue proportion are strongly regular (see Linz \cite{LinzImproved}). Hence, we are motivated to demonstrate the following:
\begin{thm}\label{thm:strreg}
If $G$ is a strongly regular graph with parameters $(n, r, e, f)$ and $n\ge 3$ then $\lambda_3(G) < \frac{n}{3}$.  
\end{thm}
We make use of the absolute bound for strongly regular graphs, which is given as Theorem 3.6.7 from Cvetković, Rowlinson and Simić \cite{SpectraTextbook} and stated as a lemma here:
\begin{lmm}[Cvetkovi\'{c} et al. \cite{SpectraTextbook}] \label{lmm:strreg}
    Let $G$ be a strongly regular graph with parameters $(n,r,e,f)$ with $n\ge 3$, such that both $G$ and $\bar{G}$ are connected. Then $G$ has spectrum $r, \lambda^{k}, \mu^{l}$ for some $r>\lambda>\mu$ and $k,l\in \mathbb{N}$ with $n=1+k+l$. Furthermore: \begin{equation*}
        n \le \min\left\{\frac{k(k+3)}{2}, \frac{l(l+3)}{2}\right\}.
    \end{equation*}
\end{lmm}
Theorem 3.1 now follows by simple casework.
\begin{proof}[Proof of Theorem \ref{thm:strreg}]
By way of contradiction, suppose $G$ is strongly regular with parameters $(n,r,e,f)$ with $n\ge 3$ and $\lambda_3(G)\ge \frac{n}{3}$. 

If $G$ is disconnected, then $f=0$. As $r(r-e-1)=(n-r-1)f$ for strongly regular graphs, either $r=0$ or $r=e+1$. In the former case, $G\cong \overline{K_n}$ has $\lambda_3(G) = 0< \frac{n}{3}$. In the latter case, we have $r+1 \mid n$ and $G\cong \frac{n}{r+1}K_{r+1}$, so $G$ has spectrum $r^{n/(r+1)},(-1)^{nr/(r+1)}$. For $\lambda_3(G) \ge \frac{n}{3}$, we require $\frac{n}{r+1}\ge 3$, but this gives $\lambda_3(G) = r \le \frac{n}{3}-1 < \frac{n}{3}$.
Similarly, if $\bar{G}$ is disconnected, then $G$ must be the $m$-fold blow-up of the complete graph $K_{n/m}$, where $m\mid n$. Thus, $G$ has spectrum $m(\frac{n}{m}-1),0^{\left(n-\frac{n}{m}\right)},(-m)^{\left(\frac{n}{m}-1\right)}$. So, $\lambda_3(G) = 0 < \frac{n}{3}$. 

Hence, we require $G$ to be connected. By Lemma \ref{lmm:strreg}, $G$ has spectrum $r,\lambda^{k},\mu^{l}$. Following the methods presented by Nikiforov \cite{NikiforovEvals}, if $A$ denotes the adjacency matrix of $G$, we have:
\begin{equation*}
    r^2 +  k\lambda^2 \le  r^2 +  k\lambda^2 + l\mu^2 = \text{Tr}(A^2) =  nr,
\end{equation*}
and therefore $\lambda^2 \le \frac{n^2}{4k}$. Hence, for $\lambda_3(G)\ge \frac{n}{3}$ we require $k\le 2$. By the second part of Lemma \ref{lmm:strreg}, we require $n\le 5$. However, the only strongly regular graph with $n\le 5$ vertices for which $G$ and $\bar{G}$ are connected is $C_5$, but $\lambda_3(C_5) = \frac{-1+\sqrt{5}}{2} < \frac{5}{3}$. 
\end{proof}

\subsection{Regular line graphs}
It is well known that the line graph of a connected graph is regular if and only if $G$ is regular or semi-regular bipartite (see Beineke and Bagga \cite{LineGraphs}). We consider these two cases separately and derive the following bounds. 
\begin{thm} If $G$ is a regular graph, then $\lambda_3(L(G)) < \frac{|V(L(G))|}{3}$.
\end{thm}
\begin{proof} Suppose that $G$ is $r$-regular on $n$ vertices with $q = \frac{1}{2}nr$ edges. Let $M,N$ denote the adjacency matrices of $G$ and $L(G)$ respectively. Let $B$ be the $n\times q$ matrix indexed by $V(G) \times E(G)$ such that:
\begin{equation*}
    B_{v,e} = \begin{cases}1 & v\in e\\ 0 & v\not\in e.\end{cases}
\end{equation*}
Then, note that $BB^T = M + rI$ and $B^TB = N + 2I$. As $BB^T$ and $B^TB$ share non-zero eigenvalues, we have $\chi_{N}(t) = (t+2)^{m-n}\chi_{M}(t-r+2)$. Thus, $\lambda_3(N)$ is either $-2$ or $\lambda_3(M)+r-2$. In any case
\begin{equation*}
    \frac{\lambda_3(L(G))}{|V(L(G))|}\le \frac{2\lambda_3(M)+2r-4}{nr} < \frac{4}{n}.
\end{equation*}
Hence, the result is immediate for $n\ge 12$. We can check the 797 line graphs of regular graphs on $n\le 11$ vertices. 
\end{proof}

\begin{thm} If $G$ is a semi-regular bipartite graph, then $\lambda_3(L(G)) \le \frac{|V(L(G))|}{3}$.
\end{thm}
\begin{proof} If $G$ is a semi-regular bipartite graph with parameters $(n_1,n_2,r_1,r_2)$ where $n_1\ge n_2$, then note that $L(G)$ is $(r_1+r_2-2)$-regular on $n_1r_1 = n_2r_2$ vertices. Hence, we are done if $\frac{n_2r_2}{3} \ge r_1 + r_2 - 2$. Namely, if $n_2 \ge 6$ then we are done as $r_1+r_2-2 < 2r_2 \le \frac{n_2 r_2}{3}$. 

If $n_2=5$, we are done unless $r_2 < \frac{3r_1-6}{2}$. As $r_1\le 5$, substituting into the previous inequality gives that $r_2\le 4$. However, we require $r_2 < \frac{3r_1-6}{2} \le \frac{3r_2-6}{2}$, which gives the contradiction that $r_2> 6$. 

If $n_2=4$, we are done unless $r_2 < 3r_1-6$. As $r_1\le 4$, we must have $r_2\le 5$. Moreover, $r_2\le 3r_2-6$ implies that $r_2\ge 3$. There are only 7 possible sets of parameters: $(5,4,4,5),(4,4,4,4),(8,4,2,4),(16,4,1,4),(4,4,3,3),(6,4,2,3),(12,4,1,3)$ and there are only 58 graphs with these parameters. We may check these manually and yield that each such line graph has third eigenvalue proportion $\le 0.202254248$.

If $n_2 = 3$, then either $G\cong 3K_{1,r_2}$, $G\cong K_{3,n_1}$ or $G$ is isomorphic to the closed vertex multiplication of $C_6$ by $[1,r_2,1,r_2,1,r_2]$. In any case, we have $\frac{\lambda_3(L(G))}{|V(L(G))|}\le \frac{1}{3}$. 

If $n_2 = 2$, we immediately have either $G\cong 2K_{1,r_2}$ or $G\cong K_{2,n_1}$. In the former case, we have $\lambda_3(L(G)) = \lambda_3(2K_{r_2}) = -1$ and in the latter, $\lambda_3(L(G)) = \lambda_3(K_{n_1} \square P_2)= 0$. Finally, if $n_2 = 1$ we have $G\cong K_{1,n_1}$ and so $\lambda_3(L(G))= - 1$. In any case, we have that $\lambda_3(L(G))=\lambda_3(K_{n_1})\le \frac{|V(L(G))|}{3}$, completing the proof.
\end{proof}

\subsection{Cayley graphs on abelian groups}

Let $\Gamma = (\Gamma, \cdot, e)$ be a group and $S$ a subset of $\Gamma$ such that $e \not \in S$, $S = S^{-1}$. The Cayley graph $G = Cay(\Gamma, S)$ is the graph with $V(G) = \Gamma$ and $s_1 \sim s_2 \iff s_1s_2^{-1} \in S$. We prove the following result and later generalise it.

\begin{thm}\label{thm:cayley}
    Suppose $G = Cay(\Gamma, S)$ where $\Gamma$ is abelian and is of order $n$. Then $\frac{\lambda_3(G)}{n} \le \frac13$. 
\end{thm}

In particular, this deals with the case of Paley graphs, which otherwise provide constructions for other extremal eigenvalue problems. We require the following three lemmas. 

\begin{lmm}\label{lmm:f}
    Let $f(n)$ denote the maximum value of $\operatorname{Re}(\sum_{\xi \in S} \xi)$ over all subsets $S$ of the $n$-th roots of unity. Let $\omega = \exp(2i\pi/n)$. Then $f(n) = \omega^{-\floor{n/4}} + \omega^{-\floor{n/4}+1} + \dots + \omega^0 + \dots + \omega^{\floor{n/4}}$. 
\end{lmm}
\begin{proof}
    Consider the $n$-th roots of unity on the complex plane. Then the ones with non-negative real part are those with argument in $[-\pi/2, \pi/2]$ hence correspond exactly to $\{\omega^{-\floor{n/4}}, \omega^{-\floor{n/4}+1}, \dots, \omega^0, \dots, \omega^{\floor{n/4}}\}$. Since $\omega^k$ and $\omega^{-k}$ are complex conjugates, $\omega^k + \omega^{-k}$ is real and the sum of this subset is real, which by construction attains maximal real part. 
\end{proof}

\begin{lmm}\label{lmm:cot}
    Let $\omega = \exp(2i\pi/n)$. Then $i\frac{\omega+1}{\omega-1} = \cot(\pi/n)$.
\end{lmm}
\begin{proof}
    We demonstrate a geometric proof. 

    \begin{figure}[H]
        \centering
        \includegraphics[scale=0.58]{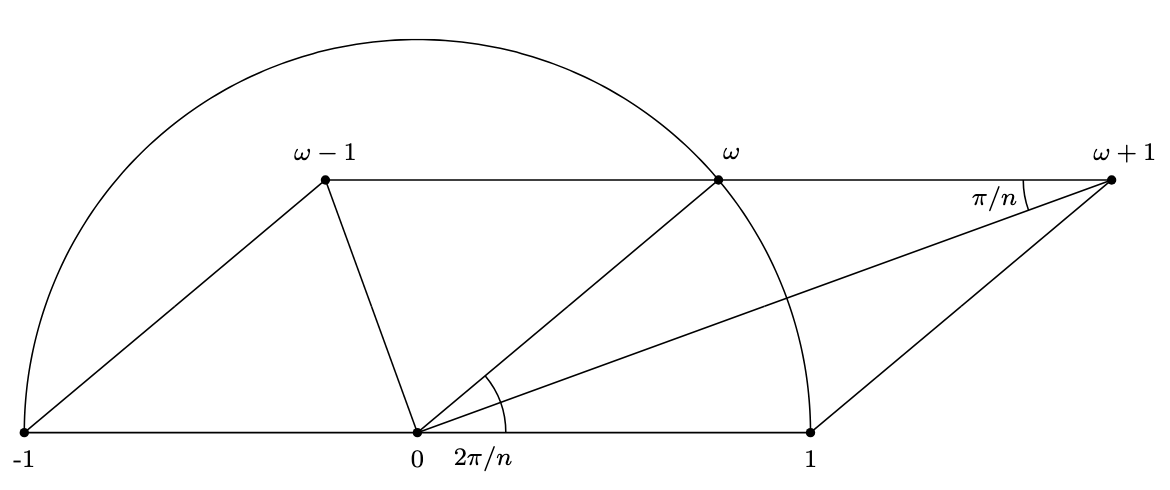}
        \caption{Two parallelograms.}
        \label{fig:parallelogram}
    \end{figure}

    Construct $\omega-1$ and $\omega+1$ in the complex plane; we have parallelograms given by $\{0, -1, \omega-1, \omega\}$ and $\{0, 1, \omega+1, \omega\}$. Thus, $\omega+1$ and $\omega-1$ are the internal and external angle bisectors of $\angle(\omega, 1)$. In particular, $\angle(\omega-1, \omega+1) = \pi/2$. Hence, $i\frac{\omega+1}{\omega-1} = \frac{|\omega+1|}{|\omega-1|} = \cot(\pi/n)$. 
\end{proof}

\begin{lmm}\label{lmm:f-growth}
    For $n \ge 3$, we have $f(n) \le n/3$ with $\lim_{n \to \infty} \frac{f(n)}{n} = \frac{1}{\pi}$.
\end{lmm}
\begin{proof}
    As a geometric series, we convert Lemma \ref{lmm:f} to $f(n) = -\omega^{\floor{n/4}} \frac{\omega^{2\floor{n/4}+1}-1}{\omega - 1}$. By computing $f(n)$ for $3 \le n \le 150$, we get $f(n) \le n/3$ with equality for $n = 3, 6$, hence we can restrict our focus to $n > 150$. 

    Consider $\delta_n = f(n) - \cot(\pi/n) = |\frac{\omega^{2\floor{n/4}+1}-1}{\omega - 1}| - |\frac{\omega^{n/2+1}-1}{\omega-1}|$. By the triangle inequality, $|\delta_n| \le \frac{|\omega^{2\floor{n/4}-n/2}-1|}{|\omega-1|}$. Note that $\alpha := 2\floor{n/4} - n/2 \in \{0, -1/2, -1, -3/2\}$. Since $n \ge 4$ and $|\sin(x)|$ is a decreasing function of $x$ in $[-\pi/2, 0]$, we have:
    \[|\omega^{2\floor{n/4}-n/2} - 1| = 2|\sin(\pi \alpha /n)| \le 2|\sin(-2\pi/n)| = |\omega^2-1|.\]
    Applying the triangle inequality again, we get $|\omega^{2\floor{n/4}-n/2} - 1| \le |\omega^2 - \omega| + |\omega - 1| = 2|\omega - 1|$. Thus, $|\delta_n| \le 2$, hence we have $\frac{f(n)}{n} \le \frac{\cot(\pi/n)+2}{n}$ and
    \[\lim_{n \to \infty} \frac{f(n)}{n} = \lim_{n \to \infty} \frac{\cot(\pi/n)}{n} = \frac{1}{\pi} \lim_{n \to \infty} \frac{\pi/n}{\tan(\pi/n)} = \frac1{\pi}.\]

    The function $\frac{\cot(\pi/x)+2}{x}$ has derivative $\frac{\pi \cot(\pi/x)^2 + \pi - x\cot(\pi/x) - 2x}{x^3}$. Note further that $\tan(x) > x$ for $x \in (0, \pi/2)$, hence $\cot(\pi/x) < x/\pi$ for $x > 2$. Now $\pi y^2 - xy = y(\pi y - x) < 0$ for $0 < y < x/\pi$, hence in particular $\pi \cot(\pi/x)^2 + \pi - x\cot(\pi/x) - 2x < \pi - 2x < 0$ for $x > 2$ and our function is decreasing. We compute for $x = 150$ the value $0.3315... < 1/3$, hence $f(n) < n/3$ for $n > 150$ as well. 
\end{proof}

We now turn to proving Theorem \ref{thm:cayley}. 

\begin{proof}[Proof of Theorem \ref{thm:cayley}]
    Since $\Gamma$ is abelian, it is equipped with $n$ irreducible characters, say $\chi_1, \chi_2, \dots, \chi_n$. Babai \cite{BabaiCayleySpectra} gave the eigenvalues of $G$ as
    \[\lambda_{\chi_i} = \sum_{s \in S} \chi_i(s), \; i = 1, 2, \dots, n.\]
    Fix a character $\chi$ and let $d$ be the largest multiplicity within the multiset $\chi(S)$. Then $\chi(\Gamma)$ is a finite subgroup of $\mathbb{C}^\times$, hence is exactly the $k$-th roots of unity for some positive integer $k$. As $d$ is the largest multiplicity, there exists $\{s_1, s_2, \dots, s_d\} \subset S$ with $\chi(s_is_1^{-1}) = 1$, hence $|\ker(\chi)| \ge d$. By the First Isomorphism Theorem, $|\ker(\chi)| = \frac{n}{|\text{Im}(\chi)|}$ which implies that $d \le \frac{n}{k}$. 

    Eigenvalues of graphs are real. We construct a partition of the multiset $\chi(S)$ into $d$ subsets of distinct elements as follows: for each element in the set $\chi(S)$, add its $i$-th copy to subset $i$. By definition of $d$, we end up with exactly $d$ subsets and each subset has real part at most $f(k)$, hence:
    \[\frac{\lambda_\chi}{n} \le \frac{df(k)}{n} \le \frac{n}{k} \cdot \frac{f(k)}{n} = \frac{f(k)}{k}.\]
    If $k \ge 3$, by Lemma \ref{lmm:f-growth} we have $\frac{\lambda_\chi}{n} \le \frac13$. If $k = 1$, we necessarily have the principal character. 
    
    Now suppose two distinct characters $\chi_1, \chi_2$ have $k = 2$. Their kernels are distinct index 2 subgroups of $\Gamma$, say $H_1, H_2$. Now there exists $x \in H_1 \setminus H_2$, hence $H_1H_2 = \Gamma$ and by the Second Isomorphism Theorem, $|H_1/H_1 \cap H_2| = |H_1H_2/H_2| = 2$, so $H_1 \cap H_2$ has index 4. 

    In particular, $n/4$ elements of $\Gamma$ have $(\chi_1(g), \chi_2(g)) = (+1, +1)$, and identically for $(+1,-1)$, $(-1,+1)$, $(-1,-1)$. Thus, $\lambda_{\chi_1} + \lambda_{\chi_2} = \sum_{s \in S} \chi_1(s) + \chi_2(s) \le n/2$ and by pigeonhole principle, one of $\lambda_{\chi_1}, \lambda_{\chi_2} \le n/4$. In particular, at most one character with $k = 2$ can have $\lambda > n/4$, so $\lambda_3(G) \le n/3$ always. 
\end{proof}

\subsection{Semi-Cayley graphs on abelian groups}

Let $\Gamma$ be a group and $R, S, T \subset \Gamma$ with $R = R^{-1}, S = S^{-1}$ and $e \not \in R \cup S$. The semi-Cayley graph $SC(\Gamma; R, S, T)$ is the graph with $V(G) = \Gamma \times \{0, 1\}$ with $(h, i) \sim (g, j)$ if and only if one of the following holds:
\begin{outline}
    \1 $i = j = 0$ and $gh^{-1} \in R$,
    \1 $i = j = 1$ and $gh^{-1} \in S$,
    \1 $i = 0, j = 1$ and $gh^{-1} \in T$.
\end{outline}
Using similar techniques to above, we prove the following result.

\begin{thm}\label{thm:semi-cayley}
    Suppose $G = SC(\Gamma; R, S, T)$ where $\Gamma$ is abelian and is of order $n$. Then $\frac{\lambda_3(G)}{2n} \le \frac13$. 
\end{thm}

As before, we prove a lemma on an analogous function closely related to $f(n)$. 

\begin{lmm}\label{lmm:g-growth}
    Let $g(n)$ denote the maximum value of $|\sum_{\xi \in S} \xi|$ over all subsets $S$ of the $n$-th roots of unity. Then $g(n) = \sqrt{f(n)^2+1}$ if $4 \mid n$, otherwise $g(n) = f(n)$. 
\end{lmm}
\begin{proof}
    Consider an arbitrary unit vector $v$. The maximal sum of projections of the $n$-th roots of unity onto $v$ occurs exactly when you take those that differ by an acute angle. Indeed, given $\{\xi^1, \xi^2, \dots, \xi^n\}$ are the roots of unity as vectors in $\mathbb{R}^2$, we have
    \[\sum_{i = 1}^n v \cdot \xi^i \le \sum_{v \cdot \xi > 0} v \cdot \xi.\]
    Note that this induces a subset of consecutive roots of unity. By Cauchy-Schwarz, 
    \[\left|v \cdot \sum_{v \cdot \xi > 0} \xi\right| \le |v|\cdot \left|\sum_{v \cdot \xi > 0} \xi\right| = \left|\sum_{v \cdot \xi > 0} \xi\right|\]
    and this sum of projections provides a lower bound on the magnitude of these roots of unity. 

    For any subset $S$ of the roots of unity, we then have $|\sum_{\xi \in S} \xi|$ is the length of its projection onto its own unit vector, which we bounded above by the magnitude of the sum of consecutive roots of unity. Thus, the maximum value of $|\sum_{\xi \in S} \xi|$ is achieved by consecutive roots, where WLOG the start is 1 as rotation doesn't affect magnitude. 

    Going clockwise around the unit circle and summing roots of unity, one constructs a regular $n$-gon. Thus, we equivalently want the maximal distance between two vertices of a regular $n$-gon with unit side length. 

    Construct its circumcircle, then the distance is increasing until the antipode then decreasing after. Hence, we require the vertices to be $\floor{n/2}$ or $\ceil{n/2}$ apart and $g(n)$ is the magnitude of the sum of either $\floor{n/2}$ or $\ceil{n/2}$ consecutive roots of unity. Considering cases on $n \pmod{4}$, we get the desired form.
\end{proof}

As in the proof of Theorem \ref{thm:cayley}, we start with an analysis of characters for which $k \ge 3$ and then consider $k = 2, 1$. 

\begin{proof}[Proof of Theorem \ref{thm:semi-cayley}]
    Since $\Gamma$ is abelian, we have $n$ irreducible characters $\chi_1, \dots, \chi_n$. Gao and Luo \cite{GaoSemiCayley} proved that the eigenvalues of $G$ corresponding to character $\chi$ are given by
    \[\lambda_\chi = \frac{\lambda_\chi^R + \lambda_\chi^S \pm \sqrt{(\lambda_\chi^R - \lambda_\chi^S)^2 + 4|\lambda_\chi^T|^2}}{2} \tag{*}\]
    where $\lambda_\chi^R = \sum_{r \in R} \chi(r)$ is the eigenvalue of $Cay(\Gamma, R)$ corresponding to $\chi$ and identically for $S$. Since $R = R^{-1}$, $S = S^{-1}$, $\lambda_\chi^R$ and $\lambda_\chi^S$ are necessarily real. We also define the analogue $\lambda_\chi^T = \sum_{t \in T} \chi(t)$ for $T$, where $\lambda_\chi^T$ need not be real. 

    Running the same multiplicity argument as in Theorem \ref{thm:cayley} and generalising Lemma \ref{lmm:f-growth} to negative sums, we get $|\lambda_\chi^R|, |\lambda_\chi^S| \le \frac{nf(k)}{k}$. Analogously, by partitioning $\chi(T)$ into subsets $T_1, \dots, T_d$ with distinct elements and noting that $|\sum_{t \in T} \chi(t)| \le |\sum_{t \in T_1} \chi(t)| + \dots + |\sum_{t \in T_d} \chi(t)|$ by the triangle inequality, Lemma \ref{lmm:g-growth} gives $|\lambda_\chi^T| \le \frac{ng(k)}{k}$. Taking the larger root in (*), we seek to bound from above $\lambda_\chi$ given these constraints. 

    Note that $\lambda_\chi$ is increasing in $|\lambda_\chi^T|$. But we also have the derivative with respect to $\lambda_\chi^R$ is $\frac12 (1 + \frac{\lambda_\chi^R - \lambda_\chi^S}{\sqrt{(\lambda_\chi^R - \lambda_\chi^S)^2 + 4|\lambda_\chi^T|^2}}) > 0$ since $|\lambda_\chi^R - \lambda_\chi^S| \le \sqrt{(\lambda_\chi^R - \lambda_\chi^S)^2 + 4|\lambda_\chi^T|^2}$. Thus, by symmetry $\lambda_\chi$ is also increasing in $\lambda_\chi^R$ and $\lambda_\chi^S$. In particular, \[\lambda_\chi \le \frac{2\frac{nf(k)}{k} + 2\frac{ng(k)}{k}}{2}\]
    and therefore
    \[\frac{\lambda_\chi}{2n} \le \frac{f(k) + g(k)}{2k}.\] 

    For $k \ge 3$ and $4 \nmid k$, this gives $\frac{\lambda_\chi}{2n} \le \frac{f(k)}{k} \le \frac13$ as before. For $k \ge 3$ and $4 \mid k$ note that we get the equality $2\floor{n/4} - n/2 = 0$ hence $\delta_n = 0$ in the proof of Lemma \ref{lmm:f-growth}. Thus $f(k) = \cot(\pi/k) \le k/\pi$. In particular, $\frac{f(k) + g(k)}{2k} \le \frac{k/\pi + \sqrt{(k/\pi)^2 + 1}}{2k} = \frac{1 + \sqrt{1 + (\pi/k)^2}}{2\pi}$. This is a decreasing function in $k$ and is less than $1/3$ for $k = 8$. By direct computation, this also holds for $k = 4$, hence for $k \ge 3$ in general we again obtain $\frac{\lambda_\chi}{2n} \le \frac13$. 

    To control $k = 2$, we first bound $\lambda_\chi$ in two cases:
    \begin{outline}[enumerate]
        \1 \textbf{Case 1:} $|\lambda_\chi^T| \le n/6$. Then note that $(\lambda_\chi^R - \lambda_\chi^S)^2 + 4|\lambda_\chi^T|^2 \le (|\lambda_\chi^R - \lambda_\chi^S| + 2|\lambda_\chi^T|)^2$ hence
        \[\lambda_\chi \le \max(\lambda_\chi^R, \lambda_\chi^S) + |\lambda_\chi^T| \le \frac{n}{2} + \frac{n}{6} = \frac{2n}{3}.\] 
        \1 \textbf{Case 2:} $|\lambda_\chi^T| \le n/3$ and $|\lambda_\chi^R - \lambda_\chi^S| \ge n/2$. Let $|\lambda_\chi^R - \lambda_\chi^S| = n/2 + t$. Then because exactly $n/2$ elements of $\Gamma$ have $\chi(g) = \pm 1$ each, $|\lambda_\chi^R|, |\lambda_\chi^S| \le n/2$ hence WLOG $\lambda_\chi^R \ge 0$ and $\lambda_\chi^S \le -t$. Then 
        \[\lambda_\chi \le \frac{n}{4} - \frac{t}{2} + \sqrt{\frac14\left(\frac{n}{2} + t\right)^2 + \frac{n^2}{9}}\]
        which has a maximum at $t = 0$, giving $\frac{2n}{3}$. 
    \end{outline}

    Consider two characters $\chi_1, \chi_2$ with $k = 2$. As in the proof of Theorem \ref{thm:cayley}, $\lambda_{\chi_1}^R + \lambda_{\chi_2}^R \le n/2$, $\lambda_{\chi_1}^S + \lambda_{\chi_2}^S \le n/2$ and individually they're bounded by $n/2$ as well. 
    
    We also claim that $|\lambda_{\chi_1}^T| + |\lambda_{\chi_2}^T| \le n/2$. Indeed, suppose there are $p$ elements of $T$ satisfying $(\chi_1(t), \chi_2(t)) = (+1, +1)$. Analogously, define $q, r, s$ for $(+1, -1), (-1, +1), (-1, -1)$ respectively. Then we get $|p + q - r - s| + |p - q + r - s| = 2u - 2v \le 2u \le n/2$ for some $u, v \subset \{p, q, r, s\}$, where $0 \le p, q, r, s \le n/4$. Thus, if we let $|\lambda_{\chi_1}^T| = y$, we have $|\lambda_{\chi_2}^T| \le n/2 - y$. 

    As before, we prove that at least one of $\chi_1, \chi_2$ gives $\lambda_\chi \le \frac{2n}{3}$ via pigeonhole. Note that the smaller eigenvalue for each $\chi$ is bounded above by $\frac{n/2 + n/2}{2} = n/2$, hence it suffices to consider only the larger one. 
    \begin{outline}
        \1[-] If $|\lambda_\chi^R - \lambda_\chi^S| \le n/2$ holds for both $\chi_1, \chi_2$, we obtain
        \[\lambda_{\chi_1} + \lambda_{\chi_2} \le \frac{n + \sqrt{(n/2)^2 + 4y^2} + \sqrt{(n/2)^2 + 4(n/2-y)^2}}{2}.\tag{*}\]
        Regarding the right-hand side of (*) as a function of $y$, we get a maximum at $y = 0$ (or $y = n/2$) where we get $\frac{3 + \sqrt{5}}{4}n < \frac{4n}{3}$. 
        \1[-] If $|\lambda_\chi^R - \lambda_\chi^S| \ge n/2$ holds for both, then we either have both satisfying $|\lambda_\chi^T| \le n/3$, or one yields $|\lambda_\chi^T| > n/3$ and the other $|\lambda_\chi^T| < n/6$, since $|\lambda_{\chi_1}^T| + |\lambda_{\chi_2}^T| \le n/2$. In the former case, by Case 2 we have $\lambda_{\chi_1}, \lambda_{\chi_2} \le \frac{2n}{3}$, and in the latter case, by Case 1 and Case 2 we again have $\lambda_{\chi_1}, \lambda_{\chi_2} \le \frac{2n}{3}$. Thus $\lambda_{\chi_1} + \lambda_{\chi_2} \le \frac{4n}{3}$.
        \1[-] If one of each, suppose WLOG $|\lambda_{\chi_1}^R - \lambda_{\chi_1}^S| = n/2 + t > n/2$. Then as in Case 2 we have say $\lambda_{\chi_1}^R \ge 0$ and $\lambda_{\chi_1}^S \le -t$, hence $\lambda_{\chi_1}^S + \lambda_{\chi_2}^S \le n/2 - t$. Combined with $\lambda_{\chi_1}^R + \lambda_{\chi_2}^R \le n/2$, we get
        \[\lambda_{\chi_1} + \lambda_{\chi_2} \le \frac{(n-t) + \sqrt{(n/2+t)^2 + 4y^2} + \sqrt{(n/2)^2 + 4(n/2-y)^2}}{2}\]
        over $0 \le t \le n/2$, $0 \le y \le n/2$. The derivative with respect to $t$ is $\frac12 (-1 + \frac{n/2+t}{\sqrt{(n/2+t)^2 + 4y^2}})$ which is always non-positive. Thus, the maximum lies at $t = 0$, reducing to (*) where $\lambda_{\chi_1} + \lambda_{\chi_2} \le \frac{3 + \sqrt{5}}{4} n < \frac{4n}{3}$. 
    \end{outline}

    The final case to eliminate is if both eigenvalues for $k = 1$, which are $\frac{|R| + |S| \pm \sqrt{(|R| - |S|)^2 + 4|T|^2}}{2}$ and one for $k = 2$ yield $\lambda_\chi > \frac{2n}{3}$. Let $\chi$ correspond to $k = 2$. Let $\frac{|R| + |S|}{2} = \frac{2n}{3} + y$, then since $\lambda_\chi > \frac{2n}{3}$, $\sqrt{(|R|-|S|)^2 + 4|T|^2} < 2y$ and $|T| < y$. Note also $|R| + |S| \le 2n$ so $y \le n/3$. If we further have $y \le n/6$ then this falls in Case 1, else if $n/6 \le y \le n/3$ with $|\lambda_\chi^R - \lambda_\chi^S| \ge n/2$ we get Case 2. Otherwise, $y > n/6$ and $|\lambda_\chi^R - \lambda_\chi^S| < n/2$. 

    Since there are $n/2$ of each of $\pm 1$ in $\chi(\Gamma)$, if $|R| = n/2+r$ then necessarily $\lambda_\chi^R \le n/2-r$ and similarly for $S$. In particular, $|R| + |S| = \frac{4n}{3} + 2y$ and therefore $\lambda_\chi^R + \lambda_\chi^S \le n - (n/3 + 2y) = \frac{2n}{3} - 2y$. Hence altogether we have
    \[\lambda_\chi \le \frac{n}{3} - y + \frac12 \sqrt{(\lambda_\chi^R - \lambda_\chi^S)^2 + 4|\lambda_\chi^T|^2} \le \frac{n}{3} - y + \frac{n}{4} + y \le \frac{2n}{3}.\qedhere\] 
\end{proof}

Gao and Luo \cite{GaoSemiCayley} note that Cayley graphs of the dihedral, $D_{2n}$, and dicyclic, $Q_{4n}$, groups are semi-Cayley graphs. Hence, we can apply the above result to these non-abelian groups.

\begin{cor}\label{cor:dihedral}
    Suppose $G = Cay(\Gamma, S)$ with $\Gamma$ of order $n$ and either dihedral or dicyclic. Then $\frac{\lambda_3(G)}{n} \le \frac13$. 
\end{cor}

\subsection{Vertex-transitive graphs}

Using the above result on Cayley graphs on abelian groups, we can generalise to vertex-transitive graphs with abelian automorphism group, as demonstrated by Lov\'{a}sz \cite{LovaszTransitive}. 

\begin{cor}\label{cor:vertex-transitive}
    Let $G$ be a vertex-transitive graph of order $n$ with abelian automorphism group $\Gamma$. Then $\frac{\lambda_3}{n} \le \frac13$. 
\end{cor}
\begin{proof}
    Lov\'{a}sz \cite{LovaszTransitive} proved that the lexicographical product of $G$ and the empty graph on $|\Gamma|/n$ vertices yields a Cayley graph on $\Gamma$, say $G'$. It has adjacency matrix $A(G') = A(G) \otimes J_{|\Gamma|/n}$, hence has spectrum $\frac{|\Gamma|}{n} \lambda_1, \dots, \frac{|\Gamma|}{n} \lambda_n, 0, \dots, 0$. If $\lambda_3 < 0$, we immediately get the conclusion and otherwise we have $\lambda_3(G') = \frac{|\Gamma|}{n} \lambda_3$. By Theorem \ref{thm:cayley}, we get $\frac{\lambda_3}{n} = \frac{|\Gamma|}{n} \frac{\lambda_3}{|\Gamma|} \le \frac13$. 
\end{proof}

\subsection{Pivalous graphs}

In an attempt to generalise the pattern of the $3K_m$ and $C_6^{[m]}$ constructions, we have found a family of graphs with surprising third eigenvalue proportion. 

\begin{dfn}\label{dfn:pivalous}
    Let $Pi_n$ denote the graph with $V(G) = \{1, 2, \dots, n\}$ and $i \sim j \iff 0 \le i - j \le \floor{\frac{n}{4}} \pmod{n}$ or $0 \le j - i \le \floor{\frac{n}{4}} \pmod{n}$. We call these the `Pivalous' graphs. 
\end{dfn}
\begin{figure}[H]
    \centering
    \includegraphics[scale=0.53]{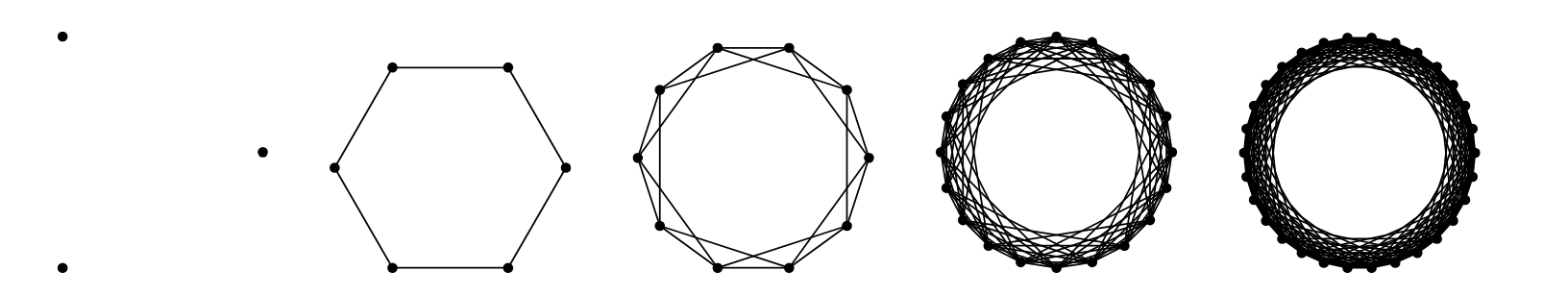}
    \caption{Examples of Pivalous graphs. From left to right: $Pi_3$, $Pi_6$, $Pi_{10}$, $Pi_{20}$, $Pi_{30}$.}
    \label{fig:pivalous-examples}
\end{figure}

This family of graphs produces circulant graphs, hence Cayley graphs on an abelian group. Note that $Pi_3 = 3K_1$ and $Pi_5, Pi_6, Pi_7$ are $C_5, C_6, C_7$ respectively. In the background of the proof of Theorem \ref{thm:cayley}, this represents the equality case. 

\begin{thm}\label{thm:pivalous-spectrum}
    The (unordered) spectrum of $Pi_n$ asymptotically grows on the order of:
    \[\frac{n}{2}, +\cot\frac{\pi}{n}, -1, -\cot\frac{3\pi}{n}, +1, +\cot\frac{5\pi}{n}, -1, \dots\]
    In particular, $\lim_{n \to \infty} \frac{\lambda_3(Pi_n)}{n} \to \frac1{\pi}$. 
\end{thm}
\begin{proof}
    Note that $Pi_n = Cay(\mathbb{Z}_n, \{-\floor{\frac{n}{4}}, \dots, \floor{\frac{n}{4}}\})$. The characters of $\mathbb{Z}_n$ are precisely $\chi_i(a) = \exp(2i\pi a/n)$. Hence, the eigenvalues of $Pi_n$ are $(\omega^t)^{-\floor{n/4}} + (\omega^t)^{-\floor{n/4}+1} + \dots + (\omega^t)^{-1} + (\omega^t)^1 + \dots + (\omega^t)^{\floor{n/4}}$ for $t = 0, 1, \dots, n-1$. For simplification, we work with $4 \mid n$ and add a self-loop at each vertex, which increases every eigenvalue by 1 and hence doesn't affect overall growth. 

    First, $\lambda_1 = n/2+1$ represents the regularity. Otherwise, we have $\lambda = (\omega^t)^{-n/4} \frac{(\omega^t)^{n/2+1} - 1}{\omega^t - 1}$. We now consider cases on $t \pmod{4}$.
    \begin{outline}
        \1[0:] $\lambda = \frac{\omega^t-1}{\omega^t-1} = +1$. 
        \1[1:] $\lambda = -i \frac{-\omega^t - 1}{\omega^t - 1} = i \frac{\omega^t + 1}{\omega^t - 1} = +\cot(\frac{t\pi}{n})$ by Lemma \ref{lmm:cot}. 
        \1[2:] $\lambda = -\frac{\omega^t-1}{\omega^t-1} = -1$. 
        \1[3:] $\lambda = i \frac{- \omega^t - 1}{\omega^t - 1} = -\cot(\frac{t\pi}{n})$.  
    \end{outline}
    Note that we get $+\cot(\frac{\pi}{n}) = -\cot(\frac{(n-1)\pi}{n})$. These are the second and third eigenvalues by size comparison to the rest, hence $\lim_{n \to \infty} \frac{\lambda_3(Pi_n)}{n} \to \frac{1}{\pi}$. 
\end{proof}

We can also consider the sequence of graphs $(Pi_n)_{n = 1}^\infty$. 

\begin{figure}[H]
    \centering
    \includegraphics[width=0.7\textwidth]{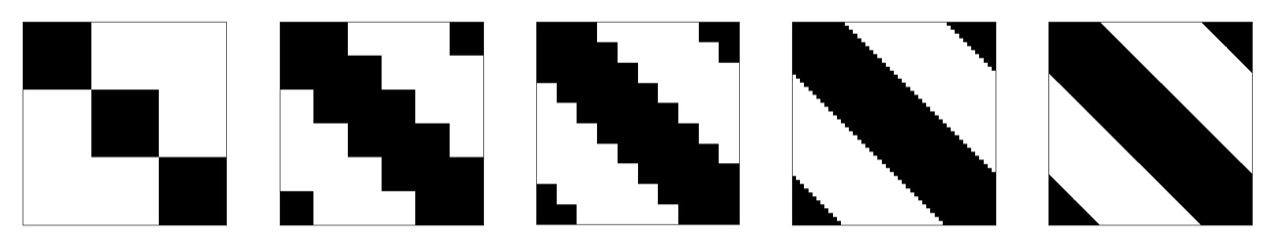}
    \caption{Visually demonstrating convergence with $Pi_3, Pi_6, Pi_{10}, Pi_{50}, Pi_{1000}$.}
    \label{fig:pivalous-limit}
\end{figure}

By observing their graphon representations in the sense of Lov\'asz and Szegedy \cite{LovaszGraphon}, we can verify that these graphs converge in cut metric to a graphon $W$ where $W(x, y) = 1$ if $0 \le x - y \le \frac14 \pmod{1}$ or $0 \le y-x \le \frac14 \pmod{1}$. Borgs et al. \cite{BorgsGraphon2} proved the limit of the spectrum is the spectrum of the limit, hence we can derive the spectrum of $W$. Alternatively, note that the eigenfunctions are $f_k(x) = \exp(2i\pi kx)$ for $k \in \mathbb{Z}$. 

\section{Extremising $\lambda_{n-1}$ over weighted graphs}\label{sec:weighted}
In this section, we consider an extension of the main problem as follows.
\begin{qun} \label{qun:weighted} Let $\mathcal{S}_n$ denote the set of $n\times n$ real symmetric matrices with entries in the interval $[0,1]$. Which matrices in $\mathcal{S}_n$ have extremal $k$-th eigenvalue; what are these extremal $k$-th eigenvalues?
\end{qun}
In other words, we consider the extension of the problem to weighted undirected graphs (possibly with loops) with weights in $[0,1]$. We remark that the function $\lambda_{k}:\mathcal{S}_n\to \mathbb{R}$ is well defined and if we equip $\mathcal{S}_n$ with the metric induced by the operator norm $|| \cdot||_{\text{opp}}$ on $n\times n$ matrices, Weyl's inequality gives:
\begin{equation*}
    |\lambda_{k}(M) - \lambda_{k}(N)| \le ||M-N||_{\text{opp}}.
\end{equation*}
Thus, $\lambda_k$ is continuous. As $\mathcal{S}_n$ is compact, $\lambda_k$ is bounded and attains its extrema, so the question is well-posed. 

Rather than considering the problem of maximising $\lambda_3$ directly, it is simpler to instead consider the problem of minimising $\lambda_{n-1}$. By Weyl's inequality we have:
\begin{equation*}
    \lambda_3(M) + \lambda_{n-1}(J-I-M) \le \lambda_{2}(J-I)=-1,
\end{equation*}
and thus,
\begin{equation*}
    \sup_{M\in \mathcal{S}_n} \lambda_3(M) \le -1 - \inf_{M\in \mathcal{S}_n} \lambda_{n-1}(M).
\end{equation*}
With this in mind, we conjecture the following:
\begin{conj} \label{conj:weighted} The function $\lambda_{n-1}: \mathcal{S}_n\to \mathbb{R}$ satisfies:
\begin{equation*}
    -\frac{n}{3} \le \inf_{M\in \mathcal{S}_n} \lambda_{n-1}(M).
\end{equation*}
\end{conj}

In this section, we make progress to answering Question \ref{qun:weighted} for $k=n-1$ and motivate Conjecture \ref{conj:weighted}. Firstly, we remark the following lemma, similar to that of Theorem 1.2 of Csikvári \cite{CsikvariConjecture}.
\begin{lmm} \label{o(n)_off_okay}
Suppose that there exists some $\kappa\in \mathbb{R}$ and $f:\mathbb{N}\to \mathbb{R}$ with $\lambda_{3}(G) \le \kappa |V(G)|+ f(|V(G)|)$ for all simple graphs $G$. Then, if  $\frac{f(n)}{n}\to 0 \text{ as } n\to \infty$, for all simple graphs $G$:
\begin{equation*}
    \lambda_{3}(G) \le \kappa |V(G)| - 1.
\end{equation*}
\end{lmm}
\begin{proof} Suppose by way of contradiction there exists some $H$ with $\lambda_3(H) = \kappa|V(H)| - 1 + \varepsilon$ with $\varepsilon>0$. Then, the closed-blowup $H^{[t]}$ of $H$ satisfies:
\begin{equation*}
   \kappa t|V(H)| + \varepsilon t - 1 = t(\lambda_3(H)+1) - 1 \le \lambda_3(H^{[t]}) \le \kappa t|V(H)| + f(t|V(H)|)
\end{equation*}
Thus, $\varepsilon \le \frac{f(t|V(H)|) + 1}{t}$ for all $t\ge 0$. But, as $t\to \infty$, $\frac{f(t|V(H)|) + 1}{t} = \frac{f(t|V(H)|) + 1}{t|V(H)|}\cdot |V(H)|\to 0$, contradicting that $\varepsilon>0$. 
\end{proof}

\subsection{Implications of the existence of an extremiser of $\lambda_{n-1}$ with $\lambda_{n}$ simple}

As previously remarked, the infimum of $\lambda_{n-1}$ is attained at some matrix $M\in \mathcal{S}_n$. Suppose that there exists such a minimiser with $\lambda_{n}(M)<\lambda_{n-1}(M)$, i.e. with $\lambda_{n}$ simple. We explore its structure should such a matrix exist.

\begin{dfn}
    Let a simple minimiser denote a matrix $M \in \mathcal{S}_n$ with $\lambda_n(M) < \lambda_{n-1}(M)$ such that $\forall M'\in \mathcal{S}_n, \lambda_{n-1}(M)\le \lambda_{n-1}(M')$. 
\end{dfn}

\begin{thm} \label{lambda_n<lambda_n-1} Suppose $M\in \mathcal{S}_n$ is a simple minimiser. Then, up to reordering rows and columns:
\begin{equation*}
    M = \begin{pmatrix} O_{P} & J_{P,N} & X \\ J_{N,P} & O_{N} & Y \\ X^T & Y^T & W
        \end{pmatrix}
\end{equation*}
where $O_m$ denotes the $m\times m$ matrix of zeros, $J_{m,l}$ denotes the $m\times l$ matrix of ones and $X,Y,W$ are $P\times Z$, $N\times Z$ and $Z\times Z$ matrices with $[0,1]$ entries respectively, with $\frac{1}{\sqrt{P}}X^T\bold{j}_P = \frac{1}{\sqrt{N}}Y^T\bold{j}_N$. If $\bold{w}$ denotes the unit eigenvector of $M$ corresponding to $\lambda_{n-1}$, then $\bold{w}^T = \begin{pmatrix} \frac{\bold{j}_P^T}{\sqrt{2P}} & -\frac{\bold{j}_N^T}{\sqrt{2N}} & \bold{0}_{Z}^T \end{pmatrix}$ and $\lambda_{n-1} = -\sqrt{PN}$.
\end{thm}
\begin{proof}
    Suppose that $M$ is a simple minimiser with eigenvalues $\lambda_n(M)<\lambda_{n-1}(M)$ and corresponding orthonormal eigenvectors $\bold{v,w}$. Suppose that $E$ is an $n\times n$ matrix such that $\exists \eta>0$ with $M+\varepsilon E \in \mathcal{S}_n$ for all $0<\varepsilon<\eta$. Then, by the min-max principle:
    \begin{alignat*}{1}
        \lambda_{n-1}(M+\varepsilon E) - \lambda_{n-1}(M) & \le \max_{\substack{\bold{x} \in \braket{\bold{v,w}} \\ ||\bold{x}|| = 1}} \left(\bold{x}^{T}(M+\varepsilon E)\bold{x} \right) - \lambda_{n-1}(M) \\
        & = \max_{\substack{\alpha, \beta \in [-1,1]\\ \alpha^2 + \beta^2 = 1}} \left(\alpha^2 (\lambda_{n}(M)-\lambda_{n-1}(M)) + \varepsilon(\alpha\bold{v} + \beta \bold{w})^{T}  E (\alpha\bold{v} + \beta \bold{w})\right)\\ & = \max_{\alpha\in [-1,1]} \left(\alpha^2 (\lambda_n(M)-\lambda_{n-1}(M) + \varepsilon(S-R)) + 2\varepsilon \alpha\sqrt{1-\alpha^2} Q + \varepsilon R\right),
    \end{alignat*}
    where $S = \bold{v}^{T}E\bold{v}$, $Q = \bold{v}^{T}E\bold{w}$ and $R = \bold{w}^{T}E\bold{w}$. Provided that $\lambda_{n-1}(M)-\lambda_{n}(M)>\varepsilon (S-R)$ (which holds for $\varepsilon$ sufficiently small), we have:
    \begin{alignat*}{1}
        \lambda_{n-1}(M+\varepsilon E) - \lambda_{n-1}(M) & = \max_{\alpha\in [-1,1]} \left(\alpha^2 (\lambda_n(M)-\lambda_{n-1}(M) + \varepsilon(S-R)) + 2\varepsilon \alpha\sqrt{1-\alpha^2} Q + \varepsilon R\right)\\ & \le  \max_{\alpha\in [-1,1]} \left(\alpha^2 (\lambda_n(M)-\lambda_{n-1}(M) + \varepsilon(S-R)) + 2\varepsilon \alpha Q + \varepsilon R\right) \\ & = \frac{\varepsilon R(\lambda_{n}(M)-\lambda_{n-1}(M) + \varepsilon(S-R)) - \varepsilon^2 Q^2}{\lambda_n(M)-\lambda_{n-1}(M) + \varepsilon(S-R)}. \tag{*}
    \end{alignat*}
    Here, we chose $\alpha$ which maximises the above quadratic.
    \begin{claim} If $E$ is a matrix such that $M+\varepsilon E \in \mathcal{S}_n$ for all $\varepsilon$ sufficiently small, then $R \ge 0$. 
    \end{claim}
    \begin{proof}
    If by way of contradiction $R < 0$ then (*) is negative for all $\varepsilon > 0$ sufficiently small. Thus, for $\varepsilon$ sufficiently small we have $\lambda_{n-1}(M+\varepsilon E)<\lambda_{n-1}(M)$, contradicting the minimality of $\lambda_{n-1}(M)$ amongst elements of ${\mathcal{S}}_n$. \end{proof}
    Let $W_{+} = \{i:w_i >0\}$, $W_{0} = \{i:w_i =0\}$ and $W_{-} = \{i:w_i<0\}$ and let $P = |W_{+}|$, $N = |W_{-}|$ and $Z = |W_{0}|$. WLOG we may reorder the vertices of the graph (i.e. rows and columns of $M$) such that $W_{+} = \{1,\dots, P\}$, $W_{-} = \{P+1,\dots, P+N\}$ and $W_{0} = \{P+N+1,\dots, P+N+Z\}$. 
    
    If by way of contradiction there exists $i,j$ with $(i,j) \in W_{+}\times W_{+} \cup W_{-}\times W_{-}$ and $M_{ij} > 0$, then consider $E$ of the form:
    \begin{equation*}
        E_{xy} = \begin{cases} -1 & \{x,y\} = \{i,j\} \\ 0 & \text{otherwise.}
        \end{cases}
    \end{equation*}
    For such $E$, $R = \bold{w}^{T}E\bold{w} = -2w_iw_j < 0$, so by the claim above, $M$ is not minimal (yielding a contradiction). Hence, if $i\ne j$ are such that $(i,j) \in W_{+}\times W_{+} \cup W_{-}\times W_{-}$, then $M_{ij} = 0$.  
    
    Similarly, if $i,j$ are such that $(i,j) \in W_{+}\times W_{-} \cup W_{-}\times W_{+}$ and $M_{ij}<1$ then considering:
    \begin{equation*}
        E_{xy} = \begin{cases} 1 & \{x,y\} = \{i,j\} \\ 0 & \text{otherwise,}
        \end{cases}
    \end{equation*}
    we again have $R = \bold{w}^{T}E\bold{w} = 2w_iw_j < 0$ giving a contradiction. Thus, if $i,j$ are such that $(i,j) \in W_{+}\times W_{-} \cup W_{-}\times W_{+}$, then $M_{ij}=1$. Hence, $M$ must have the form:
    \begin{equation*}
        M = \begin{pmatrix} O_{P} & J_{P,N} & X \\ J_{N,P} & O_{N} & Y \\ X^T & Y^T & W
        \end{pmatrix},
    \end{equation*}
    where $O_m$ denotes the $m\times m$ matrix of zeros, $J_{m,l}$ denotes the $m\times l$ matrix of ones and $X,Y,W$ are $P\times Z$, $N\times Z$ and $Z\times Z$ matrices with $[0,1]$ entries respectively. Writing $\bold{w}^T = \begin{pmatrix}\bold{w}^T_P & \bold{w}^T_N & \bold{0}_{Z}^T\end{pmatrix}$, as $\bold{w}$ is the penultimate eigenvector of $M$, we require:
    \begin{align*}
        J_{P,N}\bold{w}_{N} & = \lambda_{n-1} \bold{w}_{P} \\
        J_{N,P}\bold{w}_{P} & = \lambda_{n-1} \bold{w}_{N} \\
        X^T\bold{w}_{P} + Y^T\bold{w}_{N} & = \bold{0}_{Z}.
    \end{align*}
    Namely, the first two equations give $\bold{w}_P = x\bold{j}_P$ and $\bold{w}_N = -y\bold{j}_N$ for some $x,y>0$. Solving for $x,y$ and $\lambda_{n-1}$ gives:
    \begin{equation*}
        \lambda_{n-1} = -\sqrt{PN},\quad \bold{w}^T = \begin{pmatrix}
            \frac{\bold{j}_P^T}{\sqrt{2P}} & -\frac{\bold{j}_N^T}{\sqrt{2N}} & \bold{0}_{Z}^T
        \end{pmatrix},
    \end{equation*}
    and $\frac{X^{T}\bold{j}_{P}}{\sqrt{P}} = \frac{Y^{T}\bold{j}_{N}}{\sqrt{N}}$. 
\end{proof}

If $M$ is a simple minimiser and some $M'\in \mathcal{S}_n$ has $\lambda_{n}(M')\le \lambda_{n}(M)$ as well as $\lambda_{n-1}(M')\le\lambda_{n-1}(M)$ (so, by the minimality of $M$, $\lambda_{n-1}(M')=\lambda_{n-1}(M)$), then $M'$ is also a simple minimiser since $\lambda_n(M') \le \lambda_n(M) < \lambda_{n-1}(M) = \lambda_{n-1}(M')$. In the following, we assume the existence of such an $M$ and repeatedly construct another simple minimiser $M'$ from $M$ with nicer properties.

By the previous theorem, we have that:
    \begin{equation*}
        M = \begin{pmatrix} O_{P} & J_{P,N} & X \\ J_{N,P} & O_{N} & Y \\ X^T & Y^T & W
        \end{pmatrix},
    \end{equation*}
    with $\frac{1}{\sqrt{P}}X^T\bold{j}_P = \frac{1}{\sqrt{N}}Y^T\bold{j}_N$, $\lambda_{n-1} = -\sqrt{PN}$ and corresponding penultimate eigenvector $\bold{w}^T = \begin{pmatrix} \frac{\bold{j}_P^T}{\sqrt{2P}} & -\frac{\bold{j}_N^T}{\sqrt{2N}} & \bold{0}^T_{Z} \end{pmatrix}$. Let $\bold{v}^T = \begin{pmatrix}\bold{v}^T_P & \bold{v}^T_N & \bold{v}^T_Z\end{pmatrix}$ be a given eigenvector corresponding to $\lambda_n$. Then, without loss of generality, we may reorder the vertices $\{1,2,\dots, P\}$ so that $(\bold{v}_P)_1 \le (\bold{v}_P)_2 \le \dots \le (\bold{v}_P)_P$ and reorder the vertices $\{P+1,\dots, P+N\}$ so that $(\bold{v}_N)_1 \le (\bold{v}_N)_2 \le \dots \le (\bold{v}_N)_N$. We may also reorder the vertices $\{P+N+1,\dots P+N+Z\}$ and decompose $\bold{v}_Z^T$ as $\begin{pmatrix} \bold{v}^T_{Z_1} & \bold{v}^T_{Z_2} & \bold{v}^T_{Z_3} \end{pmatrix}$, so that $\bold{v}_{Z_1}$ has all positive entries, $\bold{v}_{Z_2}$ has all negative entries and $\bold{v}_{Z_3} = \bold{0}_{Z_3}$. Finally, decompose:
\begin{equation*}
    X = \begin{pmatrix}
        X_1 & X_2 & X_3
    \end{pmatrix} \quad Y = \begin{pmatrix}
        Y_1 & Y_2 & Y_3
    \end{pmatrix} \quad W = \begin{pmatrix}
        W^{(Z_1,Z_1)} & W^{(Z_1,Z_2)} & W^{(Z_1,Z_3)} \\ \left(W^{(Z_1,Z_2)}\right)^T & W^{(Z_2,Z_2)} & W^{(Z_2,Z_3)} \\\left(W^{(Z_1,Z_3)}\right)^T & \left(W^{(Z_2,Z_3)}\right)^T & W^{(Z_3,Z_3)}
    \end{pmatrix}
\end{equation*}
with $X_1,X_2,X_3,Y_1,Y_2,Y_3$ being $P\times Z_1, P\times Z_2, P\times Z_3$ and $N\times Z_1, N\times Z_2,N\times Z_3$ matrices with entries in $[0,1]$. 
\begin{claim} There exists a simple minimiser of the form:
\begin{equation*}
   M' = \begin{pmatrix} O_{P} & J_{P,N} & X_1 & X_2 & O_{P,Z_3} \\ J_{N,P} & O_{N} & Y_1 & Y_2 & O_{N,Z_3}\\ {X_1}^T & {Y_1}^T & W^{(Z_1,Z_1)} & W^{(Z_1,Z_2)} & O_{Z_1,Z_3} \\ {X_2}^T & {Y_2}^T & ({W^{(Z_1,Z_2)}})^T & W^{(Z_2,Z_2)}  & O_{Z_2,Z_3} \\ O_{Z_3,P} & O_{Z_3,N} & O_{Z_3,Z_1} & O_{Z_3,Z_2} & O_{Z_3}
        \end{pmatrix}.
\end{equation*}
\end{claim}
\begin{proof}
$M'$ satisfies $M'\bold{v} = M \bold{v} = \lambda_{n}\bold{v}$, $M'\bold{w} = M \bold{w} = \lambda_{n-1}\bold{w}$, hence $\lambda_{n}(M') \le \lambda_{n}(M)$ and $\lambda_{n-1}(M') \le \lambda_{n-1}(M)$.
\end{proof}
\begin{claim}
    There exists a simple minimiser of the above form with:
    \begin{equation*}
        W^{(Z_1,Z_1)} = O_{Z_1} \quad W^{(Z_1,Z_2)}= J_{Z_1,Z_2} \quad W^{(Z_2,Z_2)} = O_{Z_2}.
    \end{equation*}
\end{claim}
\begin{proof} 
Let $M$ be a simple minimiser of the form given in the previous claim, with the maximal number of 0 entries in $W^{(Z_1,Z_1)}$ and $W^{(Z_2,Z_2)}$, followed by the maximal number of 1 entries in $W^{(Z_1,Z_2)}$. Consider matrices $E$ of the form:
\begin{equation*}
    E = \begin{pmatrix}
            O_{P+N} & O_{P+N,Z_1+Z_2} & O_{P+N,Z_3} \\ O_{Z_1+Z_2,P+N} & C & O_{Z_1+Z_2,Z_3} \\ O_{Z_3,P+N} & O_{Z_3,Z_1+Z_2} & O_{Z_3}\end{pmatrix}
\end{equation*}
where $C$ is a $(Z_1+Z_2)\times (Z_1+Z_2)$ matrix such that $M+E\in \mathcal{S}_n$. We remark that $(M+E)\bold{w} = \lambda_{n-1}(M)\bold{w}$, so if $\lambda_{n}(M+E)\le \lambda_{n}(M)$ then $M+E$ is also a simple minimiser. We may bound:
\begin{equation*}
    \lambda_n(M+E) - \lambda_{n}(M) \le \bold{v}^{T}(M+E)\bold{v} - \bold{v}^{T}M\bold{v} = \bold{v}^{T}E\bold{v} = \begin{pmatrix}\bold{v}_{Z_1}^T & \bold{v}_{Z_2}^T\end{pmatrix}C\begin{pmatrix}\bold{v}_{Z_1} \\ \bold{v}_{Z_2}\end{pmatrix}.
\end{equation*}
If by way of contradiction $W^{(Z_1,Z_1)}_{ij}>0$ for some $i,j\le Z_1$, we have that taking:
\begin{equation*}
    C_{xy} = \begin{cases}
        -W^{(Z_1,Z_1)}_{ij} & \{x,y\} = \{i,j\}\\ 0 & \text{otherwise},
    \end{cases}
\end{equation*}
yields that $\begin{pmatrix}\bold{v}_{Z_1}^T&  \bold{v}_{Z_2}^{T}\end{pmatrix} C\begin{pmatrix}\bold{v}_{Z_1} \\ \bold{v}_{Z_2}\end{pmatrix} = -2W^{(Z_1,Z_1)}_{ij}(\bold{v}_{Z_1})_i(\bold{v}_{Z_1})_j<0$, so $\lambda_{n}(M+E) < \lambda_{n}(M)$. Thus, $M+E$ is another simple minimiser of $\lambda_{n-1}$ with a strictly greater number of 0 entries in $W^{(Z_1,Z_1)}$, contradicting maximality. We must thus have $W^{(Z_1,Z_1)} = O_{Z_1}$. The same logic applies to $W^{(Z_2,Z_2)}$ and we yield $W^{(Z_2,Z_2)} = O_{Z_2}$. Now, let $i\le Z_1 < j$ and suppose by way of contradiction that $W^{(Z_1,Z_2)}_{i,j-Z_1} < 1$. Then, take:
\begin{equation*}
    C_{xy} = \begin{cases}
        1-W^{(Z_1,Z_2)}_{i,j-Z_1} & \{x,y\} = \{i,j\}\\ 0 & \text{otherwise}.
    \end{cases}
\end{equation*}
As $\begin{pmatrix}\bold{v}_{Z_1}^T&  \bold{v}_{Z_2}^{T}\end{pmatrix} C\begin{pmatrix}\bold{v}_{Z_1} \\ \bold{v}_{Z_2}\end{pmatrix} = 2(1-W^{(Z_1,Z_2)}_{i,j-Z_1})(\bold{v}_{Z_1})_i(\bold{v}_{Z_2})_{j-Z_1}  < 0$, we have that $M+E$ is also a simple minimiser of $\lambda_{n-1}$. Applying this operation strictly increases the number of 1 entries in $W^{(Z_1, Z_2)}$ and preserves the entries of $W^{(Z_1, Z_1)}$ and $W^{(Z_2, Z_2)}$ -- thus we must have that $W^{(Z_1, Z_2)} = J_{Z_1,Z_2}$. 
\end{proof}
Hence, we may without loss of generality consider $M$ of the form:
\begin{equation*}
   M = \begin{pmatrix} O_{P} & J_{P,N} & X_1 & X_2 & O_{P,Z_3} \\ J_{N,P} & O_{N} & Y_1 & Y_2 & O_{N,Z_3}\\ {X_1}^T & {Y_1}^T & O_{Z_1} & J_{Z_1,Z_2} & O_{Z_1,Z_3} \\ {X_2}^T & {Y_2}^T & J_{Z_2,Z_1} & O_{Z_2}  & O_{Z_2,Z_3} \\ O_{Z_3,P} & O_{Z_3,N} & O_{Z_3,Z_1} & O_{Z_3,Z_2} & O_{Z_3}\end{pmatrix}
\end{equation*}
with $\frac{X_1^T \bold{j}_{P}}{\sqrt{P}} + \frac{Y_1^T \bold{j}_{N}}{\sqrt{N}} = 0$ and $\frac{X_2^T \bold{j}_{P}}{\sqrt{P}} + \frac{Y_2^T \bold{j}_{N}}{\sqrt{N}} = 0$. 
\begin{claim} There exists a simple minimiser of the above form with all the rows of $\begin{pmatrix} X_1^T & Y_1^T\end{pmatrix}$ equal and all the rows of $\begin{pmatrix} X_2^T & Y_2^T\end{pmatrix}$ equal.
\end{claim}
\begin{proof}
Pick $M$ a simple minimiser of the above form with the minimum number of distinct rows in $\begin{pmatrix} X_1^T & Y_1^T\end{pmatrix}$, then with the minimal number of distinct rows in $\begin{pmatrix} X_2^T & Y_2^T\end{pmatrix}$. Suppose by way of contradiction two rows of $\begin{pmatrix}X_1^T & Y_1^T\end{pmatrix}$ are not equal, say the $i$-th row and the $j$-th row. Then, consider the matrix $E$ of form:
    \begin{equation*}
        E = \begin{pmatrix}
            O_{P} & O_{P,N} & A & O_{P,Z_2} & O_{P,Z_3} \\ O_{N,P} & O_{N} & B & O_{N,Z_2} & O_{N,Z_3}\\ A^T & B^T & O_{Z_1} & O_{Z_1,Z_2} & O_{Z_1,Z_3} \\ O_{Z_2,P} & O_{Z_2,N} & O_{Z_2,Z_1} & O_{Z_2}  & O_{Z_2,Z_3} \\ O_{Z_3,P} & O_{Z_3,N} & O_{Z_3,Z_1} & O_{Z_3,Z_2} & O_{Z_3}
        \end{pmatrix}
    \end{equation*}
    where:
    \begin{equation*}
        A_{xy} = \begin{cases} (X_1)_{xi} - (X_1)_{xj} & y=j \\ 0 & \text{otherwise}\end{cases} \quad B_{xy} = \begin{cases} (Y_1)_{xi} - (Y_1)_{xj} & y=j \\ 0 & \text{otherwise.}\end{cases}
    \end{equation*}
    We remark that $(M+E)\bold{w} = -\sqrt{PN}\bold{w}$. Moreover, we have that:
    \begin{alignat*}{1}
        \lambda_{n}(M+E) - \lambda_{n}(M) & \le \bold{v}^{T} E\bold{v} = 2\sum_{x<y}  E_{xy}\bold{v}_x\bold{v}_y  \\ 
        & = 2(\bold{v}_{Z_1})_j\left(\sum_{x=1}^{P} (\bold{v}_P)_x ((X_1)_{xi} - (X_1)_{xj}) + \sum_{x=1}^{N} (\bold{v}_N)_x ((Y_1)_{xi} - (Y_1)_{xj})\right).
    \end{alignat*}
    If this is non-positive, then applying this operation to each row which is equal to the $j$-th row of $\begin{pmatrix}X_1^T&Y_1^T\end{pmatrix}$ gives another simple minimiser and strictly reduces the total number of distinct rows (contradicting maximality). Otherwise, this quantity is negative and so switching the roles of $i$ and $j$ in the above yields:
    \begin{alignat*}{1}
        \lambda_{n}(M+E) - \lambda_{n}(M) & \le \bold{v}^{T} E\bold{v} = 2\sum_{x<y}  E_{xy}\bold{v}_x\bold{v}_y  \\ 
        & = 2(\bold{v}_{Z_1})_i\left(\sum_{x=1}^{P} (\bold{v}_P)_x ((X_1)_{xj} - (X_1)_{xi}) + \sum_{x=1}^{N} (\bold{v}_N)_x ((Y_1)_{xj} - (Y_1)_{xi})\right)<0,
    \end{alignat*}
    so applying this operation to each copy of row $i$ strictly decreases $\lambda_n$ and the number of distinct rows of $\begin{pmatrix}X_1^T&Y_1^T\end{pmatrix}$, again contradicting maximality. We thus have that all rows of $\begin{pmatrix}X_1^T & Y_1^T \end{pmatrix}$ are equal. The same holds for $\begin{pmatrix}X_2^T & Y_2^T\end{pmatrix}$.
    \end{proof}
\begin{claim} 
    There exists a simple minimiser of the above form with each row of $\begin{pmatrix}X_1^T & Y_1^T\end{pmatrix}$ of the form:
    \[(\underbrace{1 \quad \dots \quad 1 \quad \alpha \quad 0 \quad \dots \quad 0}_{X_1^T} \quad \underbrace{1 \quad \dots \quad 1 \quad \beta \quad 0 \quad \dots \quad 0}_{Y_1^T}).\]
\end{claim}
\begin{proof}
Pick $M$ a simple minimiser of the above form, with the largest number of entries in $\{0,1\}$ in each row of $\begin{pmatrix}X_1^T & Y_1^T\end{pmatrix}$. Consider $E$ of the form:
\begin{equation*}
    E=\begin{pmatrix}
            O_{P} & O_{P,N} & Q \\ O_{N,P} & O_{N} & R \\ Q^T & R^T & O_{Z}\end{pmatrix}
\end{equation*}
where $Q$ and $R$ are $P\times Z$ and $N\times Z$ matrices such that $M+E\in \mathcal{S}_n$ with $\frac{Q^T\bold{j}_P}{\sqrt{P}} = \frac{R^T\bold{j}_N}{\sqrt{N}}$, the first $Z_1$ rows of $Q^T,R^T$ are equal, the next $Z_2$ are also equal and the remaining rows are 0. For such $E$, we have $(M+E)\bold{w} = -\sqrt{PN}\bold{w}$. Suppose that $p<q\le P$. Take $B=O_{N,Z}$ and:
\begin{equation*}
    Q_{xy} = \begin{cases}1 & y=p\\ -1 & y=q\\ 0 & \text{otherwise.}\end{cases}
\end{equation*}
Then, if we have $(X_1)_{ip}, (X_1)_{iq} \in (0,1)$, we obtain:
\begin{equation*}
    \lambda_{n}(M+\varepsilon E)-\lambda_{n}(M)\le \varepsilon \bold{v}^{T}E\bold{v} = 2\varepsilon \sum_{i=1}^{Z_1}(\bold{v}_{Z_1})_i((\bold{v}_{P})_p - (\bold{v}_{P})_q) \le 0
\end{equation*}
So, picking $\varepsilon=\min\{1-(X_1)_{ip},(X_1)_{iq}\}$, we have strictly increased the number of entries in $\{0,1\}$, contradicting maximality. Thus, we may suppose that there is at most one entry of each of $X_1^T, Y_1^T$ which is not $0$ nor $1$. 

Of these matrices, further pick any $M$ with the longest combined length of: the runs of 1 entries which begin at the start of each (identical) row of $X_1^T$ and $Y_1^T$ as well as the runs of 0 entries which finish at the end of each row of $X_1^T$ and $Y_1^T$. For example, the following has a combined length of $1 + 3 + 2 + 0 = 6$. 
\[(\underbrace{1 \quad 0.5 \quad 1 \quad 0 \quad 0 \quad 0}_{X_1^T} \quad \underbrace{1 \quad 1 \quad 0 \quad 0 \quad 1}_{Y_1^T}).\]

Suppose that for some $p<q\le P$ we have $(X_1)_{ip}<1$ and $(X_1)_{iq}=1$. Select the first $p$ for which this occurs. Then as above, we have that $\lambda_{n}(M+(1-(X_1)_{ip})E)-\lambda_{n}(M) = \sum_{i=1}^{Z_1} 2(1-(X_1)_{ip})(\bold{v}_{Z_1})_i((\bold{v}_P)_p - (\bold{v}_P)_q) \le 0$. The last eigenvalue does not increase, whilst the number of entries in $\{0,1\}$ stays constant and the length of the run of 1's at the start of each row of $X_1^T$ strictly increases, which gives a contradiction. We remark that this does not affect any other starting run of 1 or ending run of 0, since if $(X_1)_{ip} = 0$ then having $(X_1)_{iq} = 1$ after blocks this from connecting to the end of $X_1^T$. Thus, we have that all the 1 entries of $X_1^T$ lie in the initial run. Applying a similar argument to the run of 0 at the end of each row of $X_1^T$ and repeating for $Y_1^T$, we may conclude that each row of $\begin{pmatrix}X_1^T & Y_1^T\end{pmatrix}$ has the form:
\[(\underbrace{1 \quad \dots \quad 1 \quad \alpha \quad 0 \quad \dots \quad 0}_{X_1^T} \quad \underbrace{1 \quad \dots \quad 1 \quad \beta \quad 0 \quad \dots \quad 0}_{Y_1^T}).\qedhere\]
\end{proof}

The exact opposite holds for $\begin{pmatrix}X_2^T & Y_2^T\end{pmatrix}$ -- i.e. there exists a minimiser of the above form with each row of $\begin{pmatrix} X_2^T & Y_2^T\end{pmatrix}$ of the form:\[(\underbrace{0 \quad \dots \quad 0 \quad \alpha' \quad 1 \quad \dots \quad 1}_{X_2^T} \quad \underbrace{0 \quad \dots \quad 0 \quad \beta' \quad 1 \quad \dots \quad 1}_{Y_2^T}).\]
Hence, we may without loss of generality consider $M$ of form:
\begin{equation*}
    M = \begin{pmatrix} O_{P} & J_{P,N} & X_1 & X_2 & O_{P,Z_3} \\ J_{N,P} & O_{N} & Y_1 & Y_2 & O_{N,Z_3}\\ {X_1}^T & {Y_1}^T & O_{Z_1} & J_{Z_1,Z_2} & O_{Z_1,Z_3} \\ {X_2}^T & {Y_2}^T & J_{Z_2,Z_1} & O_{Z_2}  & O_{Z_2,Z_3} \\ O_{Z_3,P} & O_{Z_3,N} & O_{Z_3,Z_1} & O_{Z_3,Z_2} & O_{Z_3}\end{pmatrix},
\end{equation*}
with:
\begin{align*}
    X_1^T = \begin{pmatrix}J_{Z_1,Z_{1,1}} &\alpha \bold{j}_{Z_1} & O_{Z_1,Z_{1,2}}\end{pmatrix} \quad Y_1^T = \begin{pmatrix}J_{Z_1,Z_{1,3}} &\beta \bold{j}_{Z_1} & O_{Z_1,Z_{1,4}}\end{pmatrix} \\ 
    X_2^T = \begin{pmatrix}O_{Z_2,Z_{2,1}} &\alpha' \bold{j}_{Z_1} & J_{Z_2,Z_{2,2}}\end{pmatrix} \quad Y_2^T = \begin{pmatrix}O_{Z_2,Z_{2,3}} &\beta' \bold{j}_{Z_1} & J_{Z_2,Z_{2,4}}\end{pmatrix}.
\end{align*}
Consider the matrix $\tilde{M}$ given by replacing each of $\alpha,\beta,\alpha',\beta'$ with 0. It is clear that $\tilde{M}\bold{w} = -\sqrt{PN}\bold{w}$, but it may not be the case that $\bold{w}$ is the penultimate eigenvector of $\tilde{M}$. However, by Weyl's inequality, we have that:

\begin{equation*}
    |\lambda_{n-1}(M) - \lambda_{n-1}(\tilde{M})| \le \lambda_1(M-\tilde{M})\le \sqrt{\text{Tr}((M-\tilde{M})^2)} \le \sqrt{4n}.
\end{equation*}

Thus, the penultimate eigenvalue of $\tilde{M}$ differs from that of a simple minimiser by at most $2\sqrt{n}$. The advantage of working with $\tilde{M}$ is that it has a particularly simple symmetric block form, with entries in $\{0,1\}$. We remark that there are two possibilities for the overlaps between the blocks of ones and zeros of $X_1^T$ and $X_2^T$. Either the blocks of ones overlap (i.e. there exists some column index where both $X_1^{T}$ and $X_2^{T}$ have all entries given by ones) or they do not.

Accounting for both possibilities for the overlaps for $X_1^{T}$ and $X_2^{T}$ and both possibilities for the overlaps for $Y_1^T$ and $Y_2^T$, the matrix $\tilde{M}$ is thus given by a vertex multiplication (possibly with some vertices multiplied by zero) of the following $11\times 11$ matrix:
\begin{equation*}
    \setcounter{MaxMatrixCols}{20}
    T = \begin{pmatrix}
        0&0&0&0&1&1&1&1&1&0&0 \\
        0&0&0&0&1&1&1&1&1&1&0 \\
        0&0&0&0&1&1&1&1&0&0&0 \\
        0&0&0&0&1&1&1&1&0&1&0 \\
        1&1&1&1&0&0&0&0&1&0&0 \\
        1&1&1&1&0&0&0&0&1&1&0 \\
        1&1&1&1&0&0&0&0&0&0&0 \\
        1&1&1&1&0&0&0&0&0&1&0 \\
        1&1&0&0&1&1&0&0&0&1&0 \\
        0&1&0&1&0&1&0&1&1&0&0 \\
        0&0&0&0&0&0&0&0&0&0&0
    \end{pmatrix}.
\end{equation*}
Thus, we yield the following.
\begin{thm} \label{2rootn_off} If there exists a simple minimiser, then there exists some vertex multiplication of $T$, $\tilde{T}$, such that $\forall S\in \mathcal{S}_n$:
\begin{equation*}
    \lambda_{n-1}(\tilde{T}) - 2\sqrt{n} \le \lambda_{n-1}(S).
\end{equation*}
\end{thm}
We remark that the complement of $H_{a,b}$ is isomorphic to the vertex multiplication of $T$ by (for example) $[a, 0,0,b,a,0,0,b,a,b,0]$. The non-zero eigenvalues of the general vertex multiplication of $T$ are given by the roots of a sextic polynomial, which we have included in the Appendix. In light of Lemma \ref{o(n)_off_okay}, we conjecture that:
\begin{conj} For any vertex multiplication of $T$:
\begin{equation*}
    \frac{\lambda_{n-1}(\tilde{T})}{|V(\tilde{T})|} \ge -\frac{1}{3}.
\end{equation*}
\end{conj}
\subsection{Extremisers of $\lambda_{n-1}$ with $\lambda_n$ repeated} \label{sec:repfinaleig}
Throughout this section, we assume that all minimisers $M$ of $\lambda_{n-1}$ have $\lambda_{n}(M)=\lambda_{n-1}(M)$ (otherwise, we apply the results of the previous section). We give some tools to analyse the structure of such weighted graphs.
\begin{lmm} \label{negdef} Let $M$ be a minimiser of $\lambda_{n-1}$ with $\lambda_{n}(M)=\lambda_{n-1}(M)$ and $\bold{v,w}$ corresponding orthonormal eigenvectors. If there exists some symmetric $E$ such that $\bold{x}^TE\bold{x} \le 0$ for each $\bold{x}\in \braket{\bold{v,w}}$, then $\lambda_{n-1}(M+E) = \lambda_{n-1}(M)$ and if $n\ge 5$:
\begin{equation*}
    E\bold{v} = E\bold{w} = \bold{0}.
\end{equation*}
\end{lmm}
\begin{proof} Similar to Theorem \ref{lambda_n<lambda_n-1}, using min-max we have:
\begin{align*}
    0\le \lambda_{n-1}(M+E) - \lambda_{n-1}(M) & \le \sup_{\substack{\bold{x} \in \braket{\bold{v},\bold{w}}\\ ||\bold{x}|| = 1}} (\bold{x}^T(M+E)\bold{x}) - \lambda_{n-1}(M) \\ & = \sup_{\substack{\bold{x} \in \braket{\bold{v},\bold{w}}\\ ||\bold{x}|| = 1}} (\bold{x}^TE\bold{x}) \le 0.
\end{align*}
We hence require that $\lambda_{n-1}(M+E) = \lambda_{n-1}(M)$, so $M+E$ is also a minimiser of $\lambda_{n-1}$ and by the assumption at the beginning of this section, $\lambda_{n}(M+E) = \lambda_{n-1}(M+E)$. 

If $\lambda_{n-2}(M+E)=\lambda_{n-1}(M+E)$ also, then:
\begin{equation*}
    3(\lambda_{n-1}(M+E))^2 \le \text{Tr}((M+E)^2) - (\lambda_1(M+E))^2 \le \text{Tr}((M+E)^2) - \left(\frac{\text{Tr}(M+E)}{n}\right)^2 \le \frac{n^2}{4}.
\end{equation*}
Thus if $n\ge 15$, we have $\lambda_{n-1}(M+E) \ge - \frac{n}{2\sqrt{3}}> - \frac{n-2}{3}$. But, the complements of $H_{a,b}$ have $\lambda_{n-1}(\overline{H_{a,b}}) = -\frac{|V(\overline{H_{a,b}})|}{3}$, contradict the extremality of $\lambda_{n-1}(M+E)$. The authors have also verified for $5\le n\le 14$, there exists weighted graphs with penultimate eigenvalue lower than $-\frac{n}{2\sqrt{3}}$. 

Thus, $\lambda_{n-1}(M+E)< \lambda_{n-2}(M+E)$.
We also require:
\begin{equation*}
    \lambda_{n-1}(M+E) = \sup_{\substack{\bold{x} \in \braket{\bold{v},\bold{w}}\\||\bold{x}|| = 1}}\bold{x}^T (M+E)\bold{x}. \tag{*}
\end{equation*}
Let $\tilde{\bold{v}}_n, \tilde{\bold{v}}_{n-1},\dots ,\tilde{\bold{v}}_1$ be orthonormal eigenvectors of $M+E$, with corresponding eigenvalues $\tilde{\lambda}_n = \tilde{\lambda}_{n-1}, \dots \tilde{\lambda}_1$. Then, consider the vector space $V = \braket{\bold{v,w}}\cap \braket{\bold{\tilde{v}}_1,\dots \bold{\tilde{v}}_{n-1}}$. Note that $V$ has dimension at least 1, hence there exists some $\bold{z}\in V$ with $||\bold{z}||=1$. Thus, if $\bold{x} = \sum_{i=1}^{n-1} a_i \bold{\tilde{v}}_i$, then:
\begin{equation*} \bold{x}^{T}(M+E)\bold{x} =  \sum_{i=1}^{n-1} \tilde{\lambda}_{i}a_i^2 \ge \tilde{\lambda}_{n-1},
\end{equation*}
with equality iff $\bold{x} = \pm \bold{\tilde{v}}_{n-1}$. However, $\bold{x}^{T}(M+E)\bold{x} \le {\lambda}_{n-1}(M+E)$ by ($*$), so $\tilde{\bold{v}}_{n-1} \in \braket{\bold{v,w}}$. Running an identical argument but replacing $\bold{\tilde{v}}_{n-1}$ with $\bold{\tilde{v}}_n$ yields that $\tilde{\bold{v}}_{n} \in \braket{\bold{v,w}}$ and therefore $\braket{\bold{v,w}} = \braket{\tilde{\bold{v}}_{n},\tilde{\bold{v}}_{n-1}}$. Thus, $E\bold{v}=E\bold{w} = \bold{0}$ follows immediately. 
\end{proof}
We also have the following restatement of Sylvester’s criterion.
\begin{lmm} \label{SQR} Let $\bold{v,w}$ be orthonormal vectors and define $S = \bold{v}^TE\bold{v}, Q = \bold{v}^TE\bold{w}, R = \bold{w}^TE\bold{w}$. Then, $\bold{x}^TE\bold{x}\le 0$ for all $\bold{x}\in \braket{\bold{v,w}}$ iff $RS\ge Q^2$ and $R\le 0$, $S\le 0$. 
\end{lmm}
\begin{proof} We remark that:
\begin{equation*}
    \sup_{\substack{\bold{x} \in \braket{\bold{v,w}} \\ ||\bold{x}|| = 1}}\bold{x}^T E\bold{x} = \sup_{\substack{\alpha, \beta \\ \alpha^2 + \beta^2 = 1}} (S\alpha^2 + 2Q\alpha\beta + R\beta^2).
\end{equation*}
If $S\ne 0$, completing the square gives:
\begin{equation*}
    \sup_{\substack{\bold{x} \in \braket{\bold{v,w}} \\ ||\bold{x}|| = 1}}\bold{x}^T E\bold{x} = \sup_{\substack{\alpha, \beta \\ \alpha^2 + \beta^2 = 1}} \left( S\left(\alpha + \frac{Q}{S}\beta\right)^2 + \left(\frac{RS-Q^2}{S}\right)\beta^2\right),
\end{equation*}
which is less than or equal to $0$ if and only if $S< 0$ and  $RS \ge Q^2$ (implying that $R \le 0$). Identical logic applies when $R\ne 0$. 

Finally, if $R=S=0$, then:
\begin{equation*}
    \sup_{\substack{\bold{x} \in \braket{\bold{v,w}} \\ ||\bold{x}|| = 1}}\bold{x}^T E\bold{x} = \sup_{\substack{\alpha, \beta \\ \alpha^2 + \beta^2 = 1}} (2Q\alpha\beta) = |Q|,
\end{equation*}
which is less than or equal to zero iff $0=RS\ge Q^2$. 
\end{proof}
Applying Lemma \ref{negdef} and Lemma \ref{SQR} for various choices of matrix $E$ give rise to various structural properties of $M$.
\begin{lmm} \label{diagonal} If $M$ is a minimiser of $\lambda_{n-1}$, then the matrix formed by setting the diagonal entries of $M$ to 0 is also a minimiser of $\lambda_{n-1}$. 
\end{lmm}
\begin{proof} Suppose $M$ is a minimiser of $\lambda_{n-1}$ and some diagonal entry of $M$ is non-zero. Then, consider the matrix:
\begin{equation*}
    E_{xy} = \begin{cases} -M_{xx} & x=y \\ 0 & \text{otherwise.}\end{cases}
\end{equation*}
$M+E$ is precisely the matrix formed by setting the diagonal entries of $M$ to 0. Using the notation from Lemma \ref{SQR}, $S= -\sum_{i=1}^{n} M_{ii}v_i^2\le 0$, $R= -\sum_{i=1}^{n} M_{ii}w_i^2\le 0$ and $Q = -\sum_{i=1}^{n} M_{ii}v_iw_i$. Applying Cauchy-Schwarz yields $RS \ge Q^2$. The first part of Lemma \ref{negdef} thus gives $\lambda_{n-1}(M+E) = \lambda_{n-1}(M)$, so $M+E$ is also a minimiser of $\lambda_{n-1}$. 
\end{proof}
\begin{lmm} There exists a minimiser $M$ of $\lambda_{n-1}$ with at most 1 entry of each row and column not $0$ nor $1$. 
\end{lmm}
\begin{proof} Let $M$ be a minimiser of $\lambda_{n-1}$ with diagonal entries all zero and the least number of entries that are not $0$ nor $1$. Suppose by way of contradiction some row (or equivalently column) of $M$ has at least 2 entries $M_{ab}, M_{ac} \in (0,1)$ for some $a,b,c$ distinct. Consider symmetric matrices $E$ with 0 entries everywhere other than at $E_{ab},E_{ac}, E_{bc}$ (and $E_{ba}$, $E_{ca}$, $E_{cb}$). Then, we have that:
\begin{align*}
    S & = 2(E_{ab}v_av_b+E_{ac}v_av_c+E_{bc}v_bv_c) \\ 
    Q & = E_{ab}(v_aw_b+w_av_b)+E_{ac}(v_aw_c+w_av_c)+E_{bc}(v_bw_c+w_bv_c)\\ 
    R & = 2(E_{ab}w_aw_b+E_{ac}w_aw_c+E_{bc}w_bw_c).
\end{align*}
Define $[x,y]$ as $v_xw_y-v_yw_x$. Then, we have that:
\begin{align*}
    Q^2 -RS & =  \sum_{cyc} E_{ab}^2(v_aw_b-v_bw_a)^2 + 2\sum_{cyc} E_{ab}E_{bc}(v_bw_a-v_aw_b)(v_bw_c-v_cw_b).
\end{align*}
Namely, if we consider such matrices with $E_{bc}=E_{cb} = 0$, we have that $RS\ge Q^2 \iff 0 \ge (E_{ab}[a,b] - E_{ca}[c,a])^2 \iff E_{ab}[a,b] = E_{ca}[c,a]$. Moreover, $S\le 0$ iff $E_{ab}v_av_b + E_{ac}v_av_c \le 0$. Thus, provided we pick $E_{ab}\ne 0$ of the correct sign for $S\le 0$ and pick $E_{ac}$ such that $E_{ab}[a,b] = E_{ca}[c,a]$, we have that $\lambda_{n-1}(M+E)= \lambda_{n-1}(M)$. In any case, we can pick $E$ such that one of $(M+E)_{ab} \in \{0,1\}$ or $(M+E)_{ac} \in \{0,1\}$. Thus, $M+E$ is a minimiser of $\lambda_{n-1}$ with strictly fewer non-$\{0,1\}$ entries than $M$ (contradicting the hypothesis).
\end{proof}
This yields the following corollary:
\begin{cor} Regardless of the existence of a simple minimiser, there exists a symmetric $\{0,1\}$-matrix $N$ with zero entries along the diagonal such that for all $M\in \mathcal{S}_n$:
\begin{equation*}
    \lambda_{n-1}(N)\le \lambda_{n-1}(M) + 2\sqrt{n}.
\end{equation*}
Thus:
\begin{align*}
    \inf_{m\ge 2} \left\{\frac{\lambda_{m-1}(M)}{m}:M\in \mathcal{S}_m\right\} = \inf  \left\{\frac{\lambda_{|V(G)|-1}(G)}{|V(G)|}: G \text{ unweighted, } |V(G)|\ge 2\right\}.
\end{align*}
\end{cor}
We conclude this section with the following geometric constraint on the structure of minimisers $M$ of $\lambda_{n-1}$ with $\lambda_{n}(M) = \lambda_{n-1}(M)$. 
\begin{lmm} Suppose that $\lambda_{n}(M) = \lambda_{n-1}(M)=\lambda$ has corresponding orthonormal eigenvectors $\bold{v,w}$. Define the vectors $\bold{x}_1,\dots \bold{x}_n \in \mathbb{R}^2$ by:
\begin{equation*}
    \bold{x}_i = \begin{pmatrix}
        v_i \\ w_i
    \end{pmatrix}.
\end{equation*}
If there exists some $i\in \{1,\dots , n\}$, some $\lambda' \le \lambda$ and some $\alpha_1,\dots,\alpha_n \in [0,1]$ such that:
\begin{equation*}
    \lambda' \bold{x}_i = \sum_{j} \alpha_{j} \bold{x}_j,
\end{equation*}
then replacing the $i$-th row and column of $M$ with $\begin{pmatrix}\alpha_1 & \dots & \alpha_n \end{pmatrix}$ does not increase $\lambda_{n-1}$. 
\end{lmm}
\begin{proof} Consider the matrix $E$ defined by:
\begin{equation*}
    E_{xy} = \begin{cases} \alpha_{y} - M_{iy} & x=i, y\ne i \\ \alpha_{x} - M_{xi} & x\ne i,y=i \\ 0 &\text{otherwise.} \end{cases}
\end{equation*}
Then, using the notation in Lemma \ref{SQR}, we have:
\begin{align*}
    S &= 2v_i\left(\sum_{j\ne i} (\alpha_{j} - M_{ij}) v_j\right) \\
    Q &= \sum_{j\ne i} (\alpha_{j} - M_{ij}) (v_iw_j + v_jw_i) \\
    R &= 2w_i\left(\sum_{j\ne i} (\alpha_{j} - M_{ij}) w_j\right).
\end{align*}
As $\lambda' v_i = \sum_{j} \alpha_{j} v_j = \sum_{j\ne i} \alpha_{j} v_j + \alpha_i v_i$, we have that $S = 2v_i(\lambda'v_i - \alpha_i v_i) - 2\lambda v_i^2 = 2(\lambda' - \lambda - \alpha_i)v_i^2 \le 0$. Similarly, we have that $R = 2(\lambda' - \lambda - \alpha_i)w_i^2$ and:
\begin{equation*}
    Q = 2v_iw_i(\lambda' - \lambda - \alpha_i).
\end{equation*}
Thus, $RS = Q^2$. So, we have that $\lambda_{n-1}(M+E) \le \lambda_{n-1}(M)$. 
\end{proof}

\section{Concluding remarks}

With $\lambda_3 \le \frac{n}{3}$ proven for strong candidates for counter-examples, we also conjecture that $c_3 = \frac13$. Furthermore, given our work in Section \ref{sec:weighted} reducing a simple minimiser to a vertex multiplication, we are confident in the intuition that $\lambda_2 = \lambda_3$ holds for all extremal graphs, which was originally motivated by the final step in obtaining Nikiforov's bound \cite{NikiforovEvals}. Originally, we had conjectured regular to be an additional necessary property, but our discovery of $H_{a,b}$ suggests otherwise. 

An alternative perspective on $H_{a,b}$ with $a \le b$ is as an intermediate graph in the process of `morphing' a $3K_b$ into a $C_6^{[b]}$ by adding vertices to the smaller cliques in between. We also explored morphing through edges, namely turning a $3K_{2n}$ into a $C_6^{[n]}$ by splitting each $K_{2n}$ into two copies of $K_n$, pairing each copy with a copy from another $K_{2n}$, and adding edges between each pair. However, if we do this randomly with probability $p$, we obtain a graphon with spectrum 
\[\frac{p+2}{6}, \frac{1 + \sqrt{p^2 - p + 1}}{6}, \frac{1 + \sqrt{p^2 - p + 1}}{6}, \frac{1 - \sqrt{p^2 - p + 1}}{6}, \frac{1 - \sqrt{p^2 - p + 1}}{6}, \frac{-p}{6}\]
and hence third eigenvalue proportion approximately $\frac{1 + \sqrt{p^2 - p + 1}}{6}$. This only reaches $\frac13$ at $p \in \{0, 1\}$ and is worst at $p = \frac12$ which gives $\frac{2+\sqrt{3}}{12} \approx 0.311$. 

It is worth noting that our methods for Cayley graphs on abelian groups do not generalise to non-abelian, at least not easily. For example, we have by Babai \cite{BabaiCayleySpectra} that for a natural number $t$ and an irreducible character $\chi_i$ of degree $n_i$,
\[\lambda_{i,1}^t + \dots + \lambda_{i,n_i}^t = \sum_{t_1, \dots, g_t \in H} \chi_i\left(\prod_{s=1}^t g_s\right),\]
which provides a polynomial whose roots are the respective eigenvalues, but these are nonetheless difficult to extract. Regardless, due to the inherent link between sums of roots of unity and pivalous graphs, it seems plausible that the functions $f$ and $g$ would be involved again in a proof of the same result. 

Regarding Theorem \ref{2rootn_off}, although we were unable to establish explicit bounds on the penultimate root of the characteristic polynomial of the vertex multiplications of $T$, we suspect that with further simplifications or insights, such bounds could be derived. For example, if we vertex multiply $T$ by $[x_1,\dots, x_{11}]$, we see that $x_{11}$ does not affect the non-zero eigenvalues of $T$, thus we may assume that $x_{11}=0$. Note that the undirected, simple graph corresponding to $T$ has automorphism group isomorphic to $V_{4}\cong C_2\times C_2$ which simplifies matters slightly. 

Finally, while the strategy of constructing matrices $E$ for which $M+E\in \mathcal{S}_n$ and $E$ is negative semi-definite on $\braket{\bold{v,w}}$ proves effective for matrices with relatively few entries in $\{0,1\}$, this approach becomes more challenging for matrices with almost all entries in $\{0,1\}$. Nevertheless, the techniques presented provide a foothold on whether a given matrix has locally extremal penultimate eigenvalue, which we hope can be further built upon.

\section{Acknowledgements}

We pursued this research as part of the Summer Research in Mathematics (SRIM) programme organised by the Department of Pure Mathematics and Mathematical Sciences at the University of Cambridge. We were supervised by Jan Petr and are very grateful for the meet-ups, advice and insight which we found incredibly helpful. We would also like to thank the anonymous reviewers for a careful reading, valuable comments and useful suggestions.

During the programme, the first author was supported by the SRIM Bursary Fund and the second author was supported by SmartSaaS. We are thankful for their support throughout. 

\bibliography{citations}{}

\begin{thebibliography}{10}

\bibitem{BabaiCayleySpectra}
L{\'a}szl{\'o} Babai.
\newblock Spectra of {Cayley} graphs.
\newblock {\em J. Comb. Theory B}, 27:180--189, 1979.

\bibitem{LineGraphs}
Lowell~W. Beineke and Jay~S. Bagga.
\newblock {\em Spectral Properties of Line Graphs}, pages 51--59.
\newblock Springer International Publishing, Cham, 2021.

\bibitem{BorgsGraphon}
C.~Borgs, J.T. Chayes, L.~Lovász, V.T. Sós, and K.~Vesztergombi.
\newblock Convergent sequences of dense graphs {I}: Subgraph frequencies, metric properties and testing.
\newblock {\em Advances in Mathematics}, 219(6):1801--1851, 2008.

\bibitem{CsikvariConjecture}
Péter Csikvári.
\newblock On a conjecture of {V.} {Nikiforov}.
\newblock {\em Discrete Mathematics}, 309(13):4522--4526, 2009.

\bibitem{SpectraTextbook}
Dragoš Cvetković, Peter Rowlinson, and Slobodan Simić.
\newblock {\em An Introduction to the Theory of Graph Spectra}.
\newblock Cambridge University Press, New York, 2010.

\bibitem{GaoSemiCayley}
Xing Gao and Yanfeng Luo.
\newblock The spectrum of semi-{Cayley} graphs over abelian groups.
\newblock {\em Linear Algebra and its Applications}, 432(11):2974--2983, 2010.

\bibitem{LinzImproved}
William Linz.
\newblock Improved lower bounds on the extrema of eigenvalues of graphs.
\newblock {\em Graphs and Combinatorics}, 39, 07 2023.

\bibitem{LovaszTransitive}
L{\'a}szl{\'o} Lov{\'a}sz.
\newblock Spectra of graphs with transitive groups.
\newblock {\em Periodica Mathematica Hungarica}, 6:191--195, 1975.

\bibitem{LovaszGraphon}
László Lovász and Balázs Szegedy.
\newblock Szemerédi’s {Lemma} for the {Analyst}.
\newblock {\em Geometric and Functional Analysis}, 17:252--270, 04 2007.

\bibitem{NikiforovEvals}
Vladimir Nikiforov.
\newblock Extrema of graph eigenvalues.
\newblock {\em Linear Algebra and its Applications}, 482:158--190, 2015.

\bibitem{TaylorGraphs1}
D.~E. Taylor.
\newblock {\em \emph{Some topics in the theory of finite groups}}.
\newblock PhD thesis, University of Oxford, 1971.

\bibitem{TaylorGraphs2}
D.~E. Taylor.
\newblock Regular 2-graphs.
\newblock {\em Proceedings of the London Mathematical Society}, s3-35(2):257--274, 1977.

\end{thebibliography}
\bibliographystyle{plain}
\section*{Appendix}
The characteristic polynomial of the vertex multiplication $\tilde{T}$ of $T$ by $[x_1, x_2, \dots, {x}_{10}, x_{11}]$ is given by:
{\small \begin{alignat*}{2} 
\chi_{\tilde{T}}(t) = t^{11} & - && (x_{1}x_{5}+x_{1}x_{6}+x_{2}x_{5}+x_{1}x_{7}+x_{2}x_{6}+x_{3}x_{5}+x_{1}x_{8}+x_{2}x_{7}+x_{3}x_{6}+x_{4}x_{5} \\
& && +x_{1}x_{9}+x_{2}x_{8}+x_{3}x_{7}+x_{4}x_{6}+x_{2}x_{9}+x_{3}x_{8}+x_{4}x_{7}+x_{2}x_{10}+x_{4}x_{8}+x_{4}x_{10}\\ 
& && +x_{5}x_{9}+x_{6}x_{9}+x_{6}x_{10}+x_{8}x_{10}+x_{9}x_{10}) t^9 \\
&-&&2(x_{1}x_{5}x_{9}+x_{1}x_{6}x_{9}+x_{2}x_{5}x_{9}+x_{2}x_{6}x_{9}+x_{2}x_{6}x_{10}+x_{2}x_{8}x_{10}+x_{4}x_{6}x_{10} +x_{2}x_{9}x_{10} \\
& && +x_{4}x_{8}x_{10}+x_{6}x_{9}x_{10}) t^8 \\
&-&& 2(x_{1}x_{5}x_{9}+x_{1}x_{6}x_{9}+x_{2}x_{5}x_{9}+x_{2}x_{6}x_{9}+x_{2}x_{6}x_{10}+x_{2}x_{8}x_{10}+x_{4}x_{6}x_{10}+x_{2}x_{9}x_{10} \\
& && +x_{4}x_{8}x_{10}+x_{6}x_{9}x_{10}) t^8 \\ 
& + && (x_{1}x_{2}x_{5}x_{10}+x_{1}x_{3}x_{5}x_{9}+x_{1}x_{2}x_{6}x_{10}+x_{1}x_{3}x_{6}x_{9}+x_{1}x_{4}x_{5}x_{9}+x_{2}x_{3}x_{5}x_{9}+x_{1}x_{2}x_{7}x_{10} \\
& && +x_{1}x_{3}x_{7}x_{9}+x_{1}x_{4}x_{5}x_{10}+x_{1}x_{4}x_{6}x_{9}+x_{2}x_{3}x_{5}x_{10}+x_{2}x_{3}x_{6}x_{9}+x_{2}x_{4}x_{5}x_{9}+x_{1}x_{2}x_{8}x_{10}\\ 
& && +x_{1}x_{3}x_{8}x_{9}+x_{1}x_{4}x_{6}x_{10}+x_{1}x_{4}x_{7}x_{9}+x_{2}x_{3}x_{6}x_{10}+x_{2}x_{3}x_{7}x_{9}+x_{2}x_{4}x_{6}x_{9}+x_{1}x_{2}x_{9}x_{10}\\
& && +x_{1}x_{4}x_{7}x_{10}+x_{1}x_{4}x_{8}x_{9}+x_{1}x_{5}x_{6}x_{10}+x_{1}x_{5}x_{7}x_{9}+x_{2}x_{3}x_{7}x_{10}+x_{2}x_{3}x_{8}x_{9}+x_{2}x_{4}x_{7}x_{9}\\
& && +x_{3}x_{4}x_{5}x_{10}+x_{1}x_{4}x_{8}x_{10}+x_{1}x_{5}x_{8}x_{9}+x_{1}x_{6}x_{7}x_{9}+x_{2}x_{3}x_{8}x_{10}+x_{2}x_{4}x_{8}x_{9}+x_{2}x_{5}x_{6}x_{10} \\
& && +x_{2}x_{5}x_{7}x_{9}+x_{3}x_{4}x_{6}x_{10}+x_{1}x_{4}x_{9}x_{10}+x_{1}x_{5}x_{8}x_{10}+x_{1}x_{6}x_{7}x_{10}+x_{1}x_{6}x_{8}x_{9}+x_{2}x_{5}x_{8}x_{9} \\
& && +x_{2}x_{6}x_{7}x_{9}+x_{3}x_{4}x_{7}x_{10}+x_{3}x_{5}x_{6}x_{10}+x_{3}x_{5}x_{7}x_{9}+x_{1}x_{5}x_{9}x_{10}+x_{2}x_{4}x_{9}x_{10}+x_{2}x_{5}x_{8}x_{10} \\
& && +x_{2}x_{6}x_{7}x_{10}+x_{2}x_{6}x_{8}x_{9}+x_{3}x_{4}x_{8}x_{10}+x_{3}x_{5}x_{8}x_{9}+x_{3}x_{6}x_{7}x_{9}+x_{4}x_{5}x_{6}x_{10}+x_{4}x_{5}x_{7}x_{9} \\
& && +x_{1}x_{7}x_{8}x_{10}+x_{3}x_{5}x_{8}x_{10}+x_{3}x_{6}x_{7}x_{10}+x_{3}x_{6}x_{8}x_{9}+x_{4}x_{5}x_{8}x_{9}+x_{4}x_{6}x_{7}x_{9}+x_{1}x_{7}x_{9}x_{10} \\
& && +x_{2}x_{7}x_{8}x_{10}+x_{3}x_{5}x_{9}x_{10}+x_{4}x_{5}x_{8}x_{10}+x_{4}x_{6}x_{7}x_{10}+x_{4}x_{6}x_{8}x_{9}+x_{2}x_{7}x_{9}x_{10}+x_{3}x_{6}x_{9}x_{10} \\
& && +x_{3}x_{7}x_{8}x_{10}+x_{3}x_{7}x_{9}x_{10}+x_{4}x_{7}x_{8}x_{10}+x_{3}x_{8}x_{9}x_{10}+x_{4}x_{7}x_{9}x_{10}+x_{5}x_{6}x_{9}x_{10}+x_{4}x_{8}x_{9}x_{10}\\
& && +x_{5}x_{8}x_{9}x_{10}+x_{6}x_{8}x_{9}x_{10}-3x_{2}x_{6}x_{9}x_{10})t^{7} \\ 
& + && 2(x_{1}x_{2}x_{5}x_{9}x_{10}+x_{1}x_{2}x_{6}x_{9}x_{10}+x_{2}x_{3}x_{5}x_{9}x_{10}+x_{1}x_{4}x_{6}x_{9}x_{10}+x_{2}x_{3}x_{6}x_{9}x_{10}+x_{1}x_{5}x_{6}x_{9}x_{10} \\
& && +x_{2}x_{3}x_{7}x_{9}x_{10}+x_{2}x_{4}x_{6}x_{9}x_{10}+x_{2}x_{3}x_{8}x_{9}x_{10}+x_{2}x_{5}x_{6}x_{9}x_{10}+x_{1}x_{6}x_{7}x_{9}x_{10}+x_{2}x_{4}x_{8}x_{9}x_{10} \\
& && +x_{2}x_{5}x_{8}x_{9}x_{10}+x_{2}x_{6}x_{7}x_{9}x_{10}+x_{2}x_{6}x_{8}x_{9}x_{10}+x_{3}x_{6}x_{7}x_{9}x_{10}+x_{4}x_{6}x_{7}x_{9}x_{10}+x_{4}x_{6}x_{8}x_{9}x_{10} \\
& && -x_{1}x_{4}x_{7}x_{9}x_{10}-x_{3}x_{5}x_{8}x_{9}x_{10}) t^6 \\
& - && (x_{1}x_{2}x_{3}x_{5}x_{9}x_{10}+x_{1}x_{2}x_{3}x_{6}x_{9}x_{10}+x_{1}x_{2}x_{4}x_{5}x_{9}x_{10}+x_{1}x_{2}x_{3}x_{7}x_{9}x_{10}+x_{1}x_{2}x_{4}x_{6}x_{9}x_{10}\\ 
& && +x_{1}x_{3}x_{4}x_{5}x_{9}x_{10}+x_{1}x_{2}x_{3}x_{8}x_{9}x_{10}+x_{1}x_{2}x_{4}x_{7}x_{9}x_{10}+x_{1}x_{3}x_{4}x_{6}x_{9}x_{10}+x_{2}x_{3}x_{4}x_{5}x_{9}x_{10}\\
& && +x_{1}x_{2}x_{4}x_{8}x_{9}x_{10}+x_{1}x_{2}x_{5}x_{7}x_{9}x_{10}+x_{1}x_{3}x_{4}x_{7}x_{9}x_{10}+x_{1}x_{3}x_{5}x_{6}x_{9}x_{10}+x_{2}x_{3}x_{4}x_{6}x_{9}x_{10}\\
& && +x_{1}x_{2}x_{5}x_{8}x_{9}x_{10}+x_{1}x_{2}x_{6}x_{7}x_{9}x_{10}+x_{1}x_{3}x_{4}x_{8}x_{9}x_{10}+x_{1}x_{4}x_{5}x_{6}x_{9}x_{10}+x_{2}x_{3}x_{4}x_{7}x_{9}x_{10}\\
& && +x_{2}x_{3}x_{5}x_{6}x_{9}x_{10}+x_{1}x_{2}x_{6}x_{8}x_{9}x_{10}+x_{1}x_{3}x_{5}x_{8}x_{9}x_{10}+x_{1}x_{3}x_{6}x_{7}x_{9}x_{10}+x_{1}x_{4}x_{5}x_{7}x_{9}x_{10}\\
& && +x_{2}x_{3}x_{4}x_{8}x_{9}x_{10}+x_{2}x_{3}x_{5}x_{7}x_{9}x_{10}+x_{2}x_{4}x_{5}x_{6}x_{9}x_{10}+x_{1}x_{3}x_{7}x_{8}x_{9}x_{10}+x_{1}x_{4}x_{6}x_{8}x_{9}x_{10}\\
& && +x_{1}x_{5}x_{6}x_{7}x_{9}x_{10}+x_{2}x_{3}x_{6}x_{8}x_{9}x_{10}+x_{2}x_{4}x_{5}x_{8}x_{9}x_{10}+x_{2}x_{4}x_{6}x_{7}x_{9}x_{10}+x_{3}x_{4}x_{5}x_{7}x_{9}x_{10}\\
& && +x_{1}x_{4}x_{7}x_{8}x_{9}x_{10}+x_{1}x_{5}x_{6}x_{8}x_{9}x_{10}+x_{2}x_{3}x_{7}x_{8}x_{9}x_{10}+x_{2}x_{5}x_{6}x_{7}x_{9}x_{10}+x_{3}x_{4}x_{5}x_{8}x_{9}x_{10}\\
& && +x_{3}x_{4}x_{6}x_{7}x_{9}x_{10}+x_{1}x_{5}x_{7}x_{8}x_{9}x_{10}+x_{1}x_{4}x_{7}x_{8}x_{9}x_{10}+x_{2}x_{5}x_{6}x_{8}x_{9}x_{10}+x_{3}x_{4}x_{6}x_{8}x_{9}x_{10}\\
& && +x_{3}x_{5}x_{6}x_{7}x_{9}x_{10}+x_{1}x_{6}x_{7}x_{8}x_{9}x_{10}+x_{2}x_{5}x_{7}x_{8}x_{9}x_{10}+x_{3}x_{5}x_{6}x_{8}x_{9}x_{10}+x_{4}x_{5}x_{6}x_{7}x_{9}x_{10}\\
& && +x_{2}x_{6}x_{7}x_{8}x_{9}x_{10}+x_{3}x_{5}x_{7}x_{8}x_{9}x_{10}+x_{4}x_{5}x_{6}x_{8}x_{9}x_{10}+x_{3}x_{6}x_{7}x_{8}x_{9}x_{10}+x_{4}x_{5}x_{7}x_{8}x_{9}x_{10}\\
& && +x_{4}x_{6}x_{7}x_{8}x_{9}x_{10}+4x_{1}x_{4}x_{6}x_{7}x_{9}x_{10}+4x_{2}x_{3}x_{5}x_{8}x_{9}x_{10}) t^5. \\ 
\end{alignat*}}
\end{document}